# The analogue of grad-div stabilization in DG methods for incompressible flows: Limiting behavior and extension to tensor-product meshes

Mine Akbas[a], Alexander Linke[b,1], Leo G. Rebholz[c,2,*], Philipp W. Schroeder[d,3]

[a] *Department of Mathematics, Duzce University, 81620, Duzce, Turkey*
[b] *Weierstrass Institute, 10117 Berlin, Germany*
[c] *Department of Mathematical Sciences, Clemson University, Clemson, SC 29634, USA*
[d] *Institute for Numerical and Applied Mathematics, Georg-August-University Göttingen, 37083 Göttingen, Germany*

**Abstract**

Grad-div stabilization is a classical remedy in conforming mixed finite element methods for incompressible flow problems, for mitigating velocity errors that are sometimes called poor mass conservation. Such errors arise due to the relaxation of the divergence constraint in classical mixed methods, and are excited whenever the spatial discretization has to deal with comparably large and complicated pressures. In this contribution, an analogue of grad-div stabilization for Discontinuous Galerkin methods is studied. Here, the key is the penalization of the jumps of the normal velocities over facets of the triangulation, which controls the measure-valued part of the distributional divergence of the discrete velocity solution. Our contribution is twofold: first, we characterize the limit for arbitrarily large penalization parameters, which shows that the stabilized nonconforming Discontinuous Galerkin methods remain robust and accurate in this limit; second, we extend these ideas to the case of non-simplicial meshes; here, broken grad-div stabilization must be used in addition to the normal velocity jump penalization, in order to get the desired pressure robustness effect. The analysis is performed for the Stokes equations, and more complex flows and Crouzeix–Raviart elements are considered in numerical examples that also show the relevance of the theory in practical settings.

*Keywords:* Incompressible Navier–Stokes equations, mixed finite element methods, grad-div stabilization, Discontinuous Galerkin method, nonconforming finite elements

*AMS subject classifications:* 35Q30, 65M15, 65M60, 76M10

## 1 Introduction

Classical conforming and inf-sup stable mixed finite element methods for the incompressible (Navier–)Stokes equations, such as the mini [2], the Bernardi–Raugel [3] and the Taylor–Hood elements [52], relax the divergence constraint $\nabla \cdot \boldsymbol{u} = 0$, in order to construct optimally convergent spatial discretizations on regular unstructured triangulations [23, 4]. Indeed, the discrete velocity solution $\boldsymbol{u}_h$ is not divergence-free, but only *discretely divergence-free*, i.e., it holds $\nabla_h \cdot \boldsymbol{u}_h := \pi_0(\nabla \cdot \boldsymbol{u}_h) = 0$, where $\nabla_h \cdot$ denotes the discrete divergence operator and $\pi_0$ denotes the $L^2$ best approximation in the discrete pressure space.

*Corresponding author
*Email addresses:* mineakbas@duzce.edu.tr (Mine Akbas), alexander.linke@wias-berlin.de (Alexander Linke), rebholz@clemson.edu (Leo G. Rebholz), p.schroeder@math.uni-goettingen.de (Philipp W. Schroeder)
[1] ORCID: https://orcid.org/0000-0002-0165-2698
[2] supported by NSF grant DMS1522191 and U.S. Army grant 65294-MA
[3] ORCID: https://orcid.org/0000-0001-7644-4693

While relaxing the divergence constraint facilitates the construction of inf-sup stable discretizations, it was soon realized that *discretely divergence-free* is sometimes not good enough. For example, already in 1989 D. Pelletier, A. Fortin and R. Camarero titled their article [43] by the provocative question "Are FEM solutions of incompressible flows really incompressible (or how simple flows can cause headaches!)" and pointed to the problem of *poor mass conservation*. Poor mass conservation describes velocity errors that are excited when the pressure is comparably large and complicated [29]. This happens in many flow problems, including Boussinesq flows [25, 43, 17, 22, 19, 20], potential and generalized Beltrami flows [37, 36, 29], quasi-geostrophic flows [51, 14, 36], electrophoresis [44], and two-phase flows with surface tension [21, 35]. Physical problems where poor mass conservation is not strong are rather limited, and include for example pressure-driven Stokes and pressure-driven, low-Reynolds number Navier–Stokes flows through a channel with zero exterior forcing, where for the momentum balance approximately holds $-\nu\Delta\boldsymbol{u} + \nabla p \approx \mathbf{0}$, i.e., where the pressure gradient is proportional to the friction force.

Phenomenologically, poor mass conservation is often accompanied by comparably large violations of the divergence error $\|\nabla \cdot \boldsymbol{u}_h\|_{L^2}$, see for example [20, Table 1]. To address this, researchers since the late 1970's have enhanced the Navier–Stokes momentum balance of conforming mixed finite element methods by a consistent term [18, 24]

$$\gamma_{\mathrm{gd}}(\nabla \cdot \boldsymbol{u}_h, \nabla \cdot \boldsymbol{v}_h),$$

which penalizes large divergence errors, and is nowadays often called grad-div stabilization; here $\gamma_{\mathrm{gd}} \geqslant 0$ denotes the grad-div stabilization parameter. Indeed, grad-div stabilization for conforming mixed finite element methods has recently been investigated in depth, both from a theoretical and computational point of view [39, 40, 42, 20, 41, 8, 38, 28, 1]. A better understanding of grad-div stabilization was achieved, when the limit behavior for arbitrarily large stabilization parameters $\gamma_{\mathrm{gd}} \to \infty$ was investigated [8, 38, 28]; it turned out that grad-div stabilization is not so much a stabilization, but instead a kind of penalization procedure. On a fixed grid, for $\gamma_{\mathrm{gd}} \to \infty$ the grad-div stabilized discrete velocity solution $\boldsymbol{u}_h^{\gamma_{\mathrm{gd}}}$ converges to a divergence-free velocity solution $\boldsymbol{u}_h^{\infty}$ which is the solution of a *divergence-free conforming mixed finite element method* with the same discrete velocity space, but with a richer discrete pressure space. Since divergence-free, conforming mixed finite element methods are *pressure-robust* [29, 37], i.e., their velocity error does not depend on the continuous pressure, grad-div stabilized discrete velocities behave in a more robust manner against large and complicated continuous pressures. However, this theoretical understanding also revealed limitations of grad-div stabilization in that large grad-div stabilization parameters can cause classical Poisson locking phenomena, whenever the limiting divergence-free mixed method is not inf-sup stable [28]. On the other hand, on certain mesh families and for certain conforming mixed finite elements, grad-div over-stabilization can be avoided [28, Corollary 1, Case 2].

Of course, nonconforming mixed methods like the Crouzeix–Raviart finite element method [15] are as much endangered by poor mass conservation as conforming ones, since they are not pressure-robust (their velocity error depends on the continuous pressure). In [34, 33] it was recognized that pressure-robustness of a mixed method for incompressible flows does not depend on the fact that the discrete *velocity trial functions* are divergence-free, but it only depends on the *discrete velocity test functions*. They have to be divergence-free in the *weak sense* of $\boldsymbol{H}(\mathrm{div})$, in order to be orthogonal to any gradient field in the $\boldsymbol{L}^2(\Omega)$ scalar product for vector fields [29, 34]. Recall that a vector field $\boldsymbol{v} \in \boldsymbol{L}^2(\Omega)$ is said to be *weakly divergence-free*, if its *distributional divergence* [29]

$$\phi \mapsto -\int_\Omega \boldsymbol{v} \cdot \nabla\phi \, \mathrm{d}\boldsymbol{x}$$

vanishes for all $\phi \in C_0^\infty(\Omega)$, i.e., if it is orthogonal in $\boldsymbol{L}^2(\Omega)$ to all smooth gradient fields [29]. This shows that it is actually a very *strong* property for a vector field $\boldsymbol{v}$ to be *weakly divergence-free*, at least compared to being only *discretely divergence-free*. Indeed, applying the general definition of the distributional divergence to nonconforming finite element methods and Discontinuous Galerkin (DG) methods, it turns out that the distributional divergence of a velocity test function $\boldsymbol{v}_h$ is given by

$$\phi \mapsto \sum_{K \in \mathcal{T}_h} \int_K \phi \, \nabla \cdot \boldsymbol{v}_h \, \mathrm{d}\boldsymbol{x} - \sum_{F \in \mathcal{F}_h^i} \oint_F \phi(\llbracket \boldsymbol{v}_h \rrbracket \cdot \boldsymbol{n}_F) \, \mathrm{d}\boldsymbol{s},$$

since $\boldsymbol{v}_h$ is elementwise polynomial and an integration by parts can be applied. Therefore, the distributional divergence of $\boldsymbol{v}_h$ vanishes only if it holds $\nabla \cdot \boldsymbol{v}_h = 0$ elementwise for all $K \in \mathcal{T}_h$ and $\llbracket \boldsymbol{v}_h \rrbracket \cdot \boldsymbol{n}_F = 0$ for all $F \in \mathcal{F}_h^i$. Instead, if a vector field $\boldsymbol{v}_h$ is *discretely divergence-free*, this implies for usual, i.e., not pressure-robust [29], discretizations that only

$$\phi_h \mapsto \sum_{K \in \mathcal{T}_h} \int_K \phi_h \, \nabla \cdot \boldsymbol{v}_h \, \mathrm{d}\boldsymbol{x} - \sum_{F \in \mathcal{F}_h^i} \oint_F \phi_h(\llbracket \boldsymbol{v}_h \rrbracket \cdot \boldsymbol{n}_F) \, \mathrm{d}\boldsymbol{s},$$

vanishes for all $\phi_h \in Q_h$ from a *finite-dimensional* space $Q_h$ of pressure test functions. Hence, an analogue for grad-div stabilization for nonconforming methods has to penalize the elementwise defined broken divergence $\nabla_h \cdot \boldsymbol{v}_h$, *and* the facet jumps $\llbracket \boldsymbol{v}_h \rrbracket \cdot \boldsymbol{n}_F$ of the normal velocities (i.e., the mass flux). Therefore, in the following a penalization procedure with a penalization parameter $\gamma > 0$ is employed. This type of penalization for DG methods has been independently proposed in [9] and [31], and the analysis in [9] showed that it provides a benefit to velocity error analogous to what grad-div stabilization provides for conforming elements: it improves the velocity error by reducing the contribution of the pressure discretization error. In the context of the Crouzeix–Raviart finite element method, the importance of penalizing the jumps of the normal velocities was recognized by E. Burman and P. Hansbo in 2005 [7].

Herein, we make two fundamental contributions to advance the study of this type of stabilization. First, we consider the limiting behavior as the mass flux penalization parameter $\gamma \to \infty$, which is important when the viscosity is small and/or the pressure is large. We note that the limiting behavior has not yet been considered in the context of grad-div stabilization, but was considered in [13] in a similar way for constructing a 'divergence-free' HDG method applied to a 'gradient-velocity' formulation of the Stokes equations. In this limit, the discrete velocity solution $\boldsymbol{u}_h^\infty$ will be elementwise divergence-free and the jumps of the normal velocities over inner facets will vanish. Often, taking stabilization parameters large can cause over-stabilization, however we prove that the limit solution is the discrete velocity solution of other DG methods, namely of the weakly divergence-free inf-sup stable DG method proposed in [11, 12]. Hence the stabilized method is perfectly robust against over-stabilization, and moreover, this suggests that under an assumption of equal computational cost, computing the limit solution directly using the methods of [11, 12] would be advantageous over the associated stabilized DG method.

The second major contribution of this work is to extend the study of this stabilization to DG methods on tensor-product meshes. The key idea here is to use the stabilization above for the mass flux *and* an elementwise divergence penalization. In [49] an elementwise divergence penalization has been used (seemingly) for the first time. It was shown that even this incomplete stabilization method can improve poor mass conservation for both inf-sup stable $\mathbb{Q}_k^{\mathrm{dc}}/\mathbb{Q}_{k-1}^{\mathrm{dc}}$ and equal-order DG methods $\mathbb{Q}_k^{\mathrm{dc}}/\mathbb{Q}_k^{\mathrm{dc}}$ on tensor-product meshes in practical applications of non-isothermal flows, with the superscript 'dc' denoting that these elements are piecewise discontinuous. Further, it was recognized independently in the recent works [31] and [30] that in DG methods, the mass balance across interior facets has to be accounted for, *in addition to* the classical (broken) grad-div stabilization. For laminar and turbulent flows, the authors of Ref. [31] — improving earlier results [30] — use a velocity-correction time integration with equal-order $\mathbb{Q}_k^{\mathrm{dc}}/\mathbb{Q}_k^{\mathrm{dc}}$ elements on tensor-product meshes and an implementation in `deal.ii`. The additional terms they favor are slightly different from the stabilization that is considered in this work, see [31, 4.2.3 and 4.2.4]. Further, the authors do not rigorously justify their approach by numerical analysis and do not consider a possible lack of robustness against over-stabilization.

*Organization of the article:* In Section 2 we provide some notation and mathematical preliminaries to allow for a cleaner analysis to follow. Section 3 considers the case of the mass flux penalization for DG on simplices, Section 4 considers tensor product meshes, and Section 5 considers the enhancement in Crouzeix–Raviart elements. Several numerical tests of concept are given in Sections 3 and 4, and in Section 6, we consider applications of the penalization outside of the Stokes setting, to Navier–Stokes equations and to Boussinesq equations. In all numerical tests the (sometimes dramatic) improvement offered by the penalization is clear. Finally, conclusions are drawn in Section 7, and future research directions are discussed.

## 2 Stokes problem and DG setting

We consider a domain $\Omega \subset \mathbb{R}^d$, d=2,3, to be a simply connected set with smooth boundary, or a convex polygon. For $K \subseteq \Omega$ we use the standard Sobolev spaces $W^{m,p}(K)$ for scalar-valued functions with associated norms $\|\cdot\|_{W^{m,p}(K)}$ and seminorms $|\cdot|_{W^{m,p}(K)}$ for $0 \leqslant m \in \mathbb{R}$ and $p \geqslant 1$. Spaces and norms for vector- and tensor-valued functions are indicated with bold letters. We use the Lebesgue space $L^p(K) = W^{0,p}(K)$ and the Hilbert space $H^m(K) = W^{m,2}(K)$. Additionally, the closed subspaces $H_0^1(K)$ consisting of $H^1(K)$-functions with vanishing trace on $\partial K$ and the set $L_0^2(K)$ of $L^2(K)$-functions with zero mean in $K$ play an important role. The $L^2(K)$-inner product is denoted by $(\cdot,\cdot)_K$ and, if $K = \Omega$, the domain is omitted completely when no confusion can arise.

### 2.1 Continuous Stokes problem

We consider the stationary Stokes problem with no-slip boundary conditions:

$$\begin{cases} -\nu\Delta \boldsymbol{u} + \nabla p = \boldsymbol{f} & \text{in } \Omega, \\ \nabla\cdot\boldsymbol{u} = 0 & \text{in } \Omega, \\ \boldsymbol{u} = \boldsymbol{0} & \text{on } \partial\Omega. \end{cases} \tag{1}$$

With $\boldsymbol{V} = \boldsymbol{H}_0^1(\Omega)$ and $Q = L_0^2(\Omega)$, the weak formulation of (1) reads: Find $(\boldsymbol{u}, p) \in \boldsymbol{V} \times Q$ s.t., $\forall\,(\boldsymbol{v}, q) \in \boldsymbol{V} \times Q$,

$$\nu a(\boldsymbol{u},\boldsymbol{v}) + b(\boldsymbol{v},p) - b(\boldsymbol{u},q) = (\boldsymbol{f},\boldsymbol{v}). \tag{2}$$

The bilinear forms are given by

$$a(\boldsymbol{w},\boldsymbol{v}) = \int_\Omega \nabla\boldsymbol{w} : \nabla\boldsymbol{v}\,\mathrm{d}\boldsymbol{x} \quad \text{and} \quad b(\boldsymbol{w},q) = -\int_\Omega q(\nabla\cdot\boldsymbol{w})\,\mathrm{d}\boldsymbol{x}. \tag{3}$$

Weakly divergence-free velocities belong to

$$\boldsymbol{V}^{\mathrm{div}} = \{\boldsymbol{v} \in \boldsymbol{V} : b(\boldsymbol{v},q) = 0,\ \forall\, q \in Q\} = \{\boldsymbol{v} \in \boldsymbol{V} :\ \nabla\cdot\boldsymbol{v} = 0\}. \tag{4}$$

### 2.2 Discontinuous Galerkin setting

Let $\mathcal{T}_h$ be a shape-regular FE partition of $\Omega$ without hanging nodes and mesh size $h = \max_{K\in\mathcal{T}_h} h_K$, where $h_K$ denotes the diameter of the particular element $K \in \mathcal{T}_h$. Since the subsequent velocity approximation will not be $\boldsymbol{H}^1$-conforming, the broken Sobolev space is introduced as follows:

$$\boldsymbol{H}^m(\mathcal{T}_h) = \big\{\boldsymbol{w} \in \boldsymbol{L}^2(\Omega) : \boldsymbol{w}\big|_K \in \boldsymbol{H}^m(K),\ \forall\, K \in \mathcal{T}_h\big\}. \tag{5}$$

Define the broken gradient $\nabla_h : \boldsymbol{H}^1(\mathcal{T}_h) \to \boldsymbol{L}^2(\Omega)$ by

$$(\nabla_h \boldsymbol{w})\big|_K \coloneqq \nabla\big(\boldsymbol{w}\big|_K\big), \tag{6}$$

and similarly define the broken divergence. We additionally introduce the space

$$\boldsymbol{H}_0(\mathrm{div};\Omega) = \big\{\boldsymbol{w} \in \boldsymbol{L}^2(\Omega) :\ \nabla\cdot\boldsymbol{w} \in L^2(\Omega),\ \boldsymbol{w}\cdot\boldsymbol{n}\big|_{\partial\Omega} = 0\big\}. \tag{7}$$

In our context it is worth to remind the reader that the expression $\nabla\cdot\boldsymbol{w}\in L^2(\Omega)$ has the meaning that the *distributional divergence* of $\boldsymbol{w}$ can be expressed as a $L^2(\Omega)$ function, i.e., there exists $s\in L^2(\Omega)$ (called the weak divergence of $\boldsymbol{w}$) such that it holds for all $\phi\in C_0^\infty(\Omega)$ $-\int_\Omega \boldsymbol{v}\cdot\nabla\phi\,\mathrm{d}\boldsymbol{x}=\int_\Omega s\,\phi\,\mathrm{d}\boldsymbol{x}$, see [29].

The skeleton $\mathcal{F}_h$ denotes the set of all facets with $\mathcal{F}_K=\{F\in\mathcal{F}_h\colon F\subset\partial K\}$ and $h_F$ represents the diameter of each facet $F\in\mathcal{F}_h$. Note that $h_F\leqslant h_K$ holds true for all $F\in\mathcal{F}_K$ and additionally, we define $N_\partial=\max_{K\in\mathcal{T}_h}\operatorname{card}(\mathcal{F}_K)$. Moreover, $\mathcal{F}_h=\mathcal{F}_h^i\cup\mathcal{F}_h^\partial$ where $\mathcal{F}_h^i$ is the subset of interior facets and $\mathcal{F}_h^\partial$ collects all boundary facets $F\subset\partial\Omega$. To any $F\in\mathcal{F}_h$ we assign a unit normal vector $\boldsymbol{n}_F$ where, for $F\in\mathcal{F}_h^\partial$, this is the outer unit normal vector $\boldsymbol{n}$. If $F\in\mathcal{F}_h^i$, there are two adjacent elements $K^+$ and $K^-$ sharing the facet $F=\overline{\partial K^+}\cap\overline{\partial K^-}$ and $\boldsymbol{n}_F$ points in an arbitrary but fixed direction. Let $\phi$ be any piecewise smooth (scalar-, vector- or tensor-valued) function with traces from within the interior of $K^\pm$ denoted by $\phi^\pm$, respectively. Then, we define the jump $[\![\cdot]\!]_F$ and average $\{\!\{\cdot\}\!\}_F$ operator across interior facets $F\in\mathcal{F}_h^i$ by

$$[\![\phi]\!]_F=\phi^+-\phi^-\quad\text{and}\quad\{\!\{\phi\}\!\}_F=\frac{1}{2}\big(\phi^++\phi^-\big).\tag{8}$$

For boundary facets $F\in\mathcal{F}_h^\partial$ we set $[\![\phi]\!]_F=\{\!\{\phi\}\!\}_F=\phi$. These operators act componentwise for vector- and tensor-valued functions. Frequently, the subscript indicating the facet is omitted. Note that if, for example, $\boldsymbol{w}\in\boldsymbol{H}_0(\mathrm{div};\Omega)\cap\boldsymbol{H}^m(\mathcal{T}_h)$ $(m\geq 1)$, then $[\![\boldsymbol{w}]\!]\cdot\boldsymbol{n}_F=0$ for all $F\in\mathcal{F}_h^i$; cf. [16, Lemma 1.24]. This is why the stabilization of DG methods in this work is also sometimes called '$\boldsymbol{H}$(div)-stabilization'.

## 3 Mass flux penalization applied to inf-sup stable DG methods on simplicial meshes

We consider in this section analysis of DG methods on simplicial meshes, with the mass flux penalization. We note that much of Section 3.1 has essentially been done in [9] where the spatial convergence of steady Navier–Stokes with the same DG discretization and mass flux penalization is considered. We include the reduction of their analysis to the Stokes case herein because it defines the framework for our study of the limiting behavior in Section 3.4, and it gives a starting point from which to consider the case of tensor-product meshes in Section 4.

In the following, $\mathbb{P}_k(K)$ (vector-valued: $\boldsymbol{\mathbb{P}}_k(K)$) denotes the space of all polynomials on $K$ with degree less than or equal to $k$. Restricting ourselves to simplicial meshes, the global discrete spaces are

$$\boldsymbol{V}_h=\big\{\boldsymbol{v}_h\in\boldsymbol{L}^2(\Omega)\colon \boldsymbol{v}_h\big|_K\in\boldsymbol{\mathbb{P}}_k(K),\ \forall\,K\in\mathcal{T}_h\big\},\tag{9a}$$

$$Q_h=\big\{q_h\in L_0^2(\Omega)\colon q_h\big|_K\in\mathbb{P}_{k-1}(K),\ \forall\,K\in\mathcal{T}_h\big\}.\tag{9b}$$

This finite element (FE) pair is also called $\boldsymbol{\mathbb{P}}_k^{\mathrm{dc}}/\mathbb{P}_{k-1}^{\mathrm{dc}}$ and it forms a discretely inf-sup stable velocity/pressure pair; cf., for example, [46, Section 6.4]. Moreover, on simplicial meshes, we have the important property $\nabla_h\cdot\boldsymbol{\mathbb{P}}_k^{\mathrm{dc}}\subset\mathbb{P}_{k-1}^{\mathrm{dc}}$.

**Remark 3.1:** On tensor-product meshes, which we consider in the next section, the inf-sup stable pressure/velocity pair $\boldsymbol{\mathbb{Q}}_k^{\mathrm{dc}}/\mathbb{Q}_{k-1}^{\mathrm{dc}}$ is the common element choice. A key departure for such elements from the simplicial mesh framework is that $\nabla_h\cdot\boldsymbol{\mathbb{Q}}_k^{\mathrm{dc}}\not\subset\mathbb{Q}_{k-1}^{\mathrm{dc}}$, which creates an additional obstacle to handle in the analysis since here the additional penalization of the broken divergence is needed.

On $\boldsymbol{V}_h$ the following discrete trace inequality is valid; cf. [16, Remark 1.47]:

$$\forall\,\boldsymbol{v}_h\in\boldsymbol{V}_h\colon\quad \|\boldsymbol{v}_h\|_{\boldsymbol{L}^2(\partial K)}^2\leqslant C_{\mathrm{tr}}N_\partial h_K^{-1}\,\|\boldsymbol{v}_h\|_{\boldsymbol{L}^2(K)}^2,\quad\forall\,K\in\mathcal{T}_h.\tag{10}$$

A similar trace inequality holds true for the pressure space $Q_h$. The appearance below of certain traces of velocity facet values, and their normal derivatives, leads to the technical assumption for the proofs below that the involved velocities belong (at least) to $\boldsymbol{H}^{\frac{3}{2}+\varepsilon}(\mathcal{T}_h)$ for some $\varepsilon > 0$; cf. [46, Section 2.1.3]. Relaxing this assumption is possible [26], but beyond the scope of this contribution. We thus define the compound space

$$\boldsymbol{V}(h) = \boldsymbol{V}_h \oplus \left[\boldsymbol{V} \cap \boldsymbol{H}^{\frac{3}{2}+\varepsilon}(\mathcal{T}_h)\right]. \tag{11}$$

We consider the symmetric interior penalty (SIP) method with $\sigma > 0$ sufficiently large to guarantee the coercivity estimates below, and define the bilinear form

$$\begin{aligned} a_h(\boldsymbol{w}_h, \boldsymbol{v}_h) = \int_\Omega \nabla_h \boldsymbol{w}_h : \nabla_h \boldsymbol{v}_h \, \mathrm{d}\boldsymbol{x} + \sum_{F \in \mathcal{F}_h} \frac{\sigma}{h_F} \oint_F [\![\boldsymbol{w}_h]\!] \cdot [\![\boldsymbol{v}_h]\!] \, \mathrm{d}\boldsymbol{s} \\ - \sum_{F \in \mathcal{F}_h} \oint_F \{\!\{\nabla \boldsymbol{w}_h\}\!\} \boldsymbol{n}_F \cdot [\![\boldsymbol{v}_h]\!] \, \mathrm{d}\boldsymbol{s} - \sum_{F \in \mathcal{F}_h} \oint_F [\![\boldsymbol{w}_h]\!] \cdot \{\!\{\nabla \boldsymbol{v}_h\}\!\} \boldsymbol{n}_F \, \mathrm{d}\boldsymbol{s}. \end{aligned} \tag{12}$$

In fact, we assume that $\sigma \geqslant 1$, which is reasonable as this parameter usually scales with $k^2$ anyhow; cf. [16, Lemma 4.12 and Remark 1.48]. This behavior is mainly a consequence of the $k$-dependency of the discrete trace inequality (10); see, for example, [27] for more details on this topic. The natural discrete energy norm corresponding to the SIP bilinear form for $\boldsymbol{w} \in \boldsymbol{V}(h)$ is given by

$$|\!|\!|\boldsymbol{w}|\!|\!|_e^2 = \|\nabla_h \boldsymbol{w}\|_{\boldsymbol{L}^2(\Omega)}^2 + \sum_{F \in \mathcal{F}_h} \frac{\sigma}{h_F} \|[\![\boldsymbol{w}]\!]\|_{\boldsymbol{L}^2(F)}^2. \tag{13}$$

The discrete bilinear form for the pressure-velocity coupling is defined by

$$b_h(\boldsymbol{w}_h, q_h) = -\int_\Omega q_h (\nabla_h \cdot \boldsymbol{w}_h) \, \mathrm{d}\boldsymbol{x} + \sum_{F \in \mathcal{F}_h} \oint_F \{\!\{q_h\}\!\} ([\![\boldsymbol{w}_h]\!] \cdot \boldsymbol{n}_F) \, \mathrm{d}\boldsymbol{s}. \tag{14}$$

As mentioned above, the FE pair $\boldsymbol{V}_h/Q_h$ is discretely inf-sup stable. More precisely, there exists an even smaller discrete velocity space $\boldsymbol{W}_h \subset \boldsymbol{V}_h$, with $\boldsymbol{W}_h \subset \boldsymbol{H}_0(\mathrm{div};\Omega)$, such that $\boldsymbol{W}_h/Q_h$ is also inf-sup stable; cf., for example, [46, Section 6.4]. This ensures the existence of a $\beta^* > 0$, independent of $h$, such that

$$\beta^* \|q_h\|_{L^2(\Omega)} \leqslant \sup_{\boldsymbol{w}_h \in \boldsymbol{W}_h \setminus \{\boldsymbol{0}\}} \frac{b(\boldsymbol{w}_h, q_h)}{|\!|\!|\boldsymbol{w}_h|\!|\!|_e} \leqslant \sup_{\boldsymbol{v}_h \in \boldsymbol{V}_h \setminus \{\boldsymbol{0}\}} \frac{b_h(\boldsymbol{v}_h, q_h)}{|\!|\!|\boldsymbol{v}_h|\!|\!|_e}, \quad \forall\, q_h \in Q_h. \tag{15}$$

We will also utilize a stronger energy norm on $\boldsymbol{V}(h)$:

$$\forall\, \boldsymbol{w} \in \boldsymbol{V}(h)\colon \quad |\!|\!|\boldsymbol{w}|\!|\!|_{e,\sharp}^2 = |\!|\!|\boldsymbol{w}|\!|\!|_e^2 + \sum_{K \in \mathcal{T}_h} h_K \|\nabla \boldsymbol{w} \cdot \boldsymbol{n}_K\|_{\boldsymbol{L}^2(\partial K)}^2. \tag{16}$$

Then, there exists $M > 0$, independent of $h$, such that

$$\forall\, (\boldsymbol{w}, \boldsymbol{v}_h) \in \boldsymbol{V}(h) \times \boldsymbol{V}_h\colon \quad a_h(\boldsymbol{w}, \boldsymbol{v}_h) \leqslant M |\!|\!|\boldsymbol{w}|\!|\!|_{e,\sharp} |\!|\!|\boldsymbol{v}_h|\!|\!|_e. \tag{17}$$

Concerning a proof, see for example [16, Section 4.2.3] for a scalar-valued analogue. Moreover, note that the $|\!|\!|\cdot|\!|\!|_e$ and $|\!|\!|\cdot|\!|\!|_{e,\sharp}$ norms are uniformly equivalent on $\boldsymbol{V}_h$. That is, there exists a $C > 0$ such that $C|\!|\!|\boldsymbol{v}_h|\!|\!|_{e,\sharp} \leqslant |\!|\!|\boldsymbol{v}_h|\!|\!|_e \leqslant |\!|\!|\boldsymbol{v}_h|\!|\!|_{e,\sharp}$ for all $\boldsymbol{v}_h \in \boldsymbol{V}_h$; cf. [16, Lemma 4.20] (scalar-valued).

The key idea is to add a weighted form of the following mass flux penalization term to the DG formulation:

$$j_h(\boldsymbol{w}_h, \boldsymbol{v}_h) = \sum_{F \in \mathcal{F}_h} \frac{1}{h_F} \oint_F ([\![\boldsymbol{w}_h]\!] \cdot \boldsymbol{n}_F)([\![\boldsymbol{v}_h]\!] \cdot \boldsymbol{n}_F) \, \mathrm{d}\boldsymbol{s}. \tag{18}$$

We will show the remarkable positive impact this term can have, as it improves mass conservation as well as the pressure-robustness of the solution. In fact, the analytical and numerical results we obtain are similar to what is found with using grad-div stabilization in conforming methods for Stokes problems, which is why we characterize this penalization as an analogue to grad-div stabilization for nonconforming methods. This term penalizes normal jumps and therefore, roughly speaking, the difference between a fully discontinuous DG velocity and a normal-continuous $\boldsymbol{H}(\mathrm{div})$ velocity. However, we emphasize an important advantage over grad-div stabilization in conforming methods: in the proposed DG setting, over-stabilization is never possible, as it will be shown. For a discussion on the issue of over-stabilization in conforming methods with grad-div stabilization, see [28, Corollary 1].

For all $\boldsymbol{v}_h \in \boldsymbol{V}_h$ we introduce the notation

$$|\boldsymbol{v}_h|_{\mathrm{nj}}^2 = \sum_{F\in\mathcal{F}_h} \frac{1}{h_F} \oint_F (\llbracket \boldsymbol{v}_h \rrbracket \cdot \boldsymbol{n}_F)^2 \,\mathrm{d}\boldsymbol{s} = \sum_{F\in\mathcal{F}_h} \frac{1}{h_F} \left\| \llbracket \boldsymbol{v}_h \rrbracket \cdot \boldsymbol{n}_F \right\|_{L^2(F)}^2 ,$$

and note that using $\sigma \geqslant 1$, one obtains $|\boldsymbol{v}_h|_{\mathrm{nj}} \leqslant |||\boldsymbol{v}_h|||_e$ for all $\boldsymbol{v}_h \in \boldsymbol{V}(h)$.

For approximating (2), the DG method with mass flux penalization is given by: Find $(\boldsymbol{u}_h, p_h) \in \boldsymbol{V}_h \times Q_h$ satisfying $\forall\, (\boldsymbol{v}_h, q_h) \in \boldsymbol{V}_h \times Q_h$,

$$\nu a_h(\boldsymbol{u}_h, \boldsymbol{v}_h) + b_h(\boldsymbol{v}_h, p_h) + \gamma j_h(\boldsymbol{u}_h, \boldsymbol{v}_h) = (\boldsymbol{f}, \boldsymbol{v}_h), \tag{19}$$

$$-b_h(\boldsymbol{u}_h, q_h) = 0, \tag{20}$$

where $\gamma \geqslant 0$ is the penalization parameter.

Discretely divergence-free DG velocities belong to

$$\boldsymbol{V}_h^{\mathrm{div}} = \{\boldsymbol{v}_h \in \boldsymbol{V}_h \colon b_h(\boldsymbol{v}_h, q_h) = 0,\ \forall\, q_h \in Q_h\},$$

and thanks to (15), an equivalent and pressure-free formulation of (19)–(20) can be expressed as:

$$\text{Find } \boldsymbol{u}_h \in \boldsymbol{V}_h^{\mathrm{div}} \text{ s.t. } \nu a_h(\boldsymbol{u}_h, \boldsymbol{v}_h) + \gamma j_h(\boldsymbol{u}_h, \boldsymbol{v}_h) = (\boldsymbol{f}, \boldsymbol{v}_h), \quad \forall\, \boldsymbol{v}_h \in \boldsymbol{V}_h^{\mathrm{div}}.$$

### 3.1 Energy estimate

Provided $\sigma \geqslant 1$ is sufficiently large, making use of the discrete coercivity property of the SIP bilinear form $a_h$, cf., for example, [46, Lemma 6.6] or [16, Section 6.1.2.1], we also easily obtain discrete coercivity with a constant $C_\sigma > 0$, independent of $h$, in the following sense:

$$\forall\, \boldsymbol{v}_h \in \boldsymbol{V}_h \colon \quad a_h(\boldsymbol{v}_h, \boldsymbol{v}_h) + \gamma j_h(\boldsymbol{v}_h, \boldsymbol{v}_h) \geqslant C_\sigma |||\boldsymbol{v}_h|||_e^2 + \gamma |\boldsymbol{v}_h|_{\mathrm{nj}}^2 \geqslant C_\sigma |||\boldsymbol{v}_h|||_e^2. \tag{21}$$

The well-posedness of the formulation thus follows using this discrete coercivity property and discrete inf-sup stability. The following energy estimate follows immediately as well.

**Lemma 3.2 (*Energy estimate*)**

*Let $\boldsymbol{f} \in \boldsymbol{L}^2(\Omega)$ and assume that $\sigma > 0$ is sufficiently large to guarantee discrete coercivity. Then, with a constant $C > 0$, one obtains the following estimate for the FEM solution $(\boldsymbol{u}_h, p_h)$ to (19)–(20):*

$$\frac{\nu C_\sigma}{2} |||\boldsymbol{u}_h|||_e^2 + \gamma |\boldsymbol{u}_h|_{\mathrm{nj}}^2 \leqslant C\nu^{-1} \left\| \boldsymbol{f} \right\|_{\boldsymbol{L}^2(\Omega)}^2 , \tag{22a}$$

$$\left\| \nabla_h \cdot \boldsymbol{u}_h \right\|_{L^2(\Omega)}^2 \leqslant C |\boldsymbol{u}_h|_{\mathrm{nj}}^2 \leqslant C\gamma^{-1}\nu^{-1} \left\| \boldsymbol{f} \right\|_{\boldsymbol{L}^2(\Omega)}^2 , \tag{22b}$$

$$\left\| p_h \right\|_{L^2(\Omega)}^2 \leqslant C \left\| \boldsymbol{f} \right\|_{\boldsymbol{L}^2(\Omega)}^2 . \tag{22c}$$

**Proof:** Testing with $(\boldsymbol{v}_h, q_h) = (\boldsymbol{u}_h, p_h)$ in (19)–(20), together with coercivity from (21) and Cauchy–Schwarz leads to

$$\nu C_\sigma \|\boldsymbol{u}_h\|_e^2 + \gamma |\boldsymbol{u}_h|_{\mathrm{nj}}^2 \leqslant \|\boldsymbol{f}\|_{\boldsymbol{L}^2(\Omega)} \|\boldsymbol{u}_h\|_{\boldsymbol{L}^2(\Omega)}. \tag{23}$$

Further estimating the right-hand side requires a DG analogue of the Poincaré–Friedrichs (PF) inequality; cf., for example, [16, Corollary 5.4]. Then, Young's inequality can be invoked to obtain

$$\|\boldsymbol{f}\|_{\boldsymbol{L}^2(\Omega)} \|\boldsymbol{u}_h\|_{\boldsymbol{L}^2(\Omega)} \leqslant \frac{1}{2}\frac{C_{\mathrm{PF}}}{\nu C_\sigma} \|\boldsymbol{f}\|_{\boldsymbol{L}^2(\Omega)}^2 + \frac{\nu C_\sigma}{2} \|\boldsymbol{u}_h\|_e^2. \tag{24}$$

Reordering shows the first bound. For an estimate of the divergence of the DG solution, relying on $\nabla_h \cdot \mathbb{P}_k^{\mathrm{dc}} \subset \mathbb{P}_{k-1}^{\mathrm{dc}}$ makes it possible to choose $q_h = \nabla_h \cdot \boldsymbol{u}_h$ in (20). After an additional application of the Cauchy–Schwarz inequality and inserting $\frac{h_F}{\gamma}\frac{\gamma}{h_F} = 1$, we have

$$\begin{aligned}
\|\nabla_h \cdot \boldsymbol{u}_h\|_{L^2(\Omega)}^2 &= \sum_{F\in\mathcal{F}_h} \oint_F \{\!\{ \nabla \cdot \boldsymbol{u}_h \}\!\} ([\![\boldsymbol{u}_h]\!] \cdot \boldsymbol{n}_F)\, \mathrm{d}\boldsymbol{s} \\
&\leqslant \left( \sum_{F\in\mathcal{F}_h} \frac{h_F}{\gamma} \|\{\!\{ \nabla \cdot \boldsymbol{u}_h \}\!\}\|_{L^2(F)}^2 \right)^{1/2} \left( \sum_{F\in\mathcal{F}_h} \frac{\gamma}{h_F} \|[\![\boldsymbol{u}_h]\!] \cdot \boldsymbol{n}_F\|_{L^2(F)}^2 \right)^{1/2}.
\end{aligned} \tag{25}$$

For any $q_h \in Q_h$, the pressure analogue of the discrete trace inequality (10) yields

$$\begin{aligned}
\sum_{F\in\mathcal{F}_h} \oint_F \left|\{\!\{ q_h \}\!\}\right|^2 \mathrm{d}\boldsymbol{s} &\leqslant \sum_{F\in\mathcal{F}_h} \left[ \|q_h^+\|_{L^2(F)}^2 + \|q_h^-\|_{L^2(F)}^2 \right] \\
&\leqslant \sum_{K\in\mathcal{T}_h} \|q_h\|_{L^2(\partial K)}^2 \\
&\leqslant \sum_{K\in\mathcal{T}_h} C_{\mathrm{tr}} N_\partial h_K^{-1} \|q_h\|_{L^2(K)}^2.
\end{aligned}$$

Now, we insert $\{\!\{ q_h \}\!\} = \{\!\{ \nabla \cdot \boldsymbol{u}_h \}\!\}$ in this estimate and multiply the integrand by $h_F/\gamma$. Thus, with a generic constant $C > 0$, using $h_F \leqslant h_K$ for all $F \in \mathcal{F}_K$, we obtain the following estimate for the first term on the right-hand side of (25):

$$\left( \sum_{F\in\mathcal{F}_h} \frac{h_F}{\gamma} \|\{\!\{ \nabla \cdot \boldsymbol{u}_h \}\!\}\|_{L^2(F)}^2 \right)^{1/2} \leqslant C\gamma^{-1/2} \|\nabla_h \cdot \boldsymbol{u}_h\|_{L^2(\Omega)}. \tag{26}$$

Simply inserting (26) into (25) shows the first estimate in (22b). However, in order to further estimate the second term on the right-hand side of (25), we use the energy estimate (22a) for the velocity; i.e.,

$$\left( \sum_{F\in\mathcal{F}_h} \frac{\gamma}{h_F} \|[\![\boldsymbol{u}_h]\!] \cdot \boldsymbol{n}_F\|_{L^2(F)}^2 \right)^{1/2} = \gamma^{1/2} |\boldsymbol{u}_h|_{\mathrm{nj}} \leqslant C\nu^{-1/2} \|\boldsymbol{f}\|_{\boldsymbol{L}^2(\Omega)}. \tag{27}$$

Using (27) and (26) in (25) shows the second bound in (22b).

For the pressure bound, we invoke the discrete inf-sup condition (15) in $\boldsymbol{W}_h/Q_h$ and obtain from (19) that

$$\beta^* \|p_h\|_{L^2(\Omega)} \leqslant \sup_{\boldsymbol{w}_h \in \boldsymbol{W}_h \setminus \{\boldsymbol{0}\}} \frac{b_h(\boldsymbol{w}_h, p_h)}{\|\boldsymbol{w}_h\|_e} = \sup_{\boldsymbol{w}_h \in \boldsymbol{W}_h \setminus \{\boldsymbol{0}\}} \left[ \frac{(\boldsymbol{f}, \boldsymbol{w}_h)}{\|\boldsymbol{w}_h\|_e} - \nu \frac{a_h(\boldsymbol{u}_h, \boldsymbol{w}_h)}{\|\boldsymbol{w}_h\|_e} \right]. \tag{28}$$

Here, the fact that $j_h(\boldsymbol{u}_h, \boldsymbol{w}_h) = 0$ since $\boldsymbol{w}_h \in \boldsymbol{W}_h \subset \boldsymbol{H}_0(\mathrm{div};\Omega)$ has been used. A further estimation of the right-hand side uses Cauchy–Schwarz and Poincaré–Friedrichs for the term involving $\boldsymbol{f}$ together with boundedness of $a_h$ (17) and uniform equivalence of $|||\cdot|||_e$ and $|||\cdot|||_{e,\sharp}$ on $\boldsymbol{V}_h$ for the viscous term. Thus, the supremum vanishes thereby leading to

$$\beta^* \|p_h\|_{L^2(\Omega)} \leqslant C[\|\boldsymbol{f}\|_{\boldsymbol{L}^2} + \nu|||\boldsymbol{u}_h|||_e]. \tag{29}$$

Inserting the energy estimate (22a) for the velocity concludes the proof. ∎

### 3.2 Error estimate

For the error analysis that follows, it is important that the following property holds. A proof is straightforward and based on the consistency of both the DG method and $\boldsymbol{H}(\mathrm{div})$-stabilization.

**Corollary 3.3 (*Galerkin orthogonality*)**

*Let $(\boldsymbol{u}, p) \in \boldsymbol{V} \times Q$ solve (2) and $(\boldsymbol{u}_h, p_h) \in \boldsymbol{V}_h \times Q_h$ solve (19)–(20). If additionally $\boldsymbol{u} \in \boldsymbol{H}^{\frac{3}{2}+\varepsilon}(\mathcal{T}_h)$ and $p \in \boldsymbol{H}^{\frac{1}{2}+\varepsilon}(\mathcal{T}_h)$ for $\varepsilon > 0$, then, for all $(\boldsymbol{v}_h, q_h) \in \boldsymbol{V}_h \times Q_h$:*

$$\nu a_h(\boldsymbol{u} - \boldsymbol{u}_h, \boldsymbol{v}_h) + b_h(\boldsymbol{v}_h, p - p_h) - b_h(\boldsymbol{u} - \boldsymbol{u}_h, q_h) + \gamma j_h(\boldsymbol{u} - \boldsymbol{u}_h, \boldsymbol{v}_h) = 0. \tag{30}$$

Under the assumptions of the previous corollary we decompose the error as

$$\begin{aligned} \boldsymbol{u} - \boldsymbol{u}_h &= (\boldsymbol{u} - \boldsymbol{\pi}_h\boldsymbol{u}) - (\boldsymbol{u}_h - \boldsymbol{\pi}_h\boldsymbol{u}) = \boldsymbol{\eta}^{\boldsymbol{u}} - \boldsymbol{e}_h^{\boldsymbol{u}}, \\ p - p_h &= (p - \pi_0 p) - (p_h - \pi_0 p) = \eta^p - e_h^p, \end{aligned}$$

where $(\boldsymbol{\pi}_h, \pi_0) \colon \boldsymbol{V} \times Q \to \boldsymbol{V}_h \times Q_h$ represents appropriate approximation operators and we refer to $(\boldsymbol{\eta}^{\boldsymbol{u}}, \eta^p)$ and $(\boldsymbol{e}_h^{\boldsymbol{u}}, e_h^p)$ as approximation error and discretization error, respectively. For the pressure, $\pi_0$ simply denotes the local $L^2$-projection onto $Q_h$. For the approximation operator for the velocity, we require $\boldsymbol{\pi}_h \colon \boldsymbol{V} \to \boldsymbol{W}_h \subset \boldsymbol{V}_h$, and recall from above, where the inf-sup stability condition (15) has been introduced, that $\boldsymbol{W}_h \subset \boldsymbol{H}_0(\mathrm{div};\Omega)$. It is well-known that functions belonging to $\boldsymbol{H}_0(\mathrm{div};\Omega)$ have continuous normal components across interior facets. Specifically, we define the operator $\boldsymbol{\pi}_h$ to be the Brezzi–Douglas–Marini (BDM) interpolation operator of order $k$; cf. [4, Sections 2.3 and 2.5]. A very important property of the $\mathrm{BDM}_k$ interpolator $\boldsymbol{\pi}_h$ is the following commuting diagram property:

$$\forall\, \boldsymbol{w} \in \boldsymbol{H}_0(\mathrm{div};\Omega)\colon \qquad \nabla\cdot(\boldsymbol{\pi}_h \boldsymbol{w}) = \pi_0(\nabla\cdot\boldsymbol{w}). \tag{31}$$

Note that for the Stokes velocity $\boldsymbol{u}$, we have that $\nabla\cdot\boldsymbol{u} = 0$ holds weakly, which implies that $\nabla\cdot(\boldsymbol{\pi}_h\boldsymbol{u}) = 0$ also holds weakly. This property is important for the following theorem.

**Theorem 3.4 (*Error estimate*)**

*Let $\boldsymbol{f} \in \boldsymbol{L}^2(\Omega)$ and assume that $\sigma > 0$ is sufficiently large to guarantee discrete coercivity. Under the assumptions of the previous corollary, the following error estimates hold true:*

$$|||\boldsymbol{u} - \boldsymbol{u}_h|||_e \leqslant C\Big[|||\boldsymbol{\eta}^{\boldsymbol{u}}|||_{e,\sharp} + \gamma^{-1/2}\nu^{-1/2}\,\|\eta^p\|_{L^2(\Omega)}\Big] \tag{32a}$$

$$\|\pi_0 p - p_h\|_{L^2(\Omega)} \leqslant C\nu|||\boldsymbol{u} - \boldsymbol{u}_h|||_{e,\sharp} \tag{32b}$$

$$\|p - p_h\|_{L^2(\Omega)} \leqslant C\nu|||\boldsymbol{\eta}^{\boldsymbol{u}}|||_{e,\sharp} + \left[C\sqrt{\frac{\nu}{\gamma}} + 1\right]\|\eta^p\|_{L^2(\Omega)} \tag{32c}$$

**Remark 3.5:** Due to the fact that the BDM interpolator has optimal approximation properties, one obtains the standard convergence rate of $h^k$ whenever the exact solution $(\boldsymbol{u}, p)$ is smooth enough. Thus, the overall spatial convergence rate, as $h \to 0$, of the stabilized method remains the same as the unstabilized method, see e.g. [46]. Moreover, the accuracy of the discretization is robust and optimal with respect to the limit $\gamma \to \infty$, i.e., over-stabilization is not possible.

**Remark 3.6:** Although the mass flux penalization does not alter the spatial convergence rate, it can dramatically lower the velocity error for certain flows. In particular, in the unstabilized DG method, the coefficient of the pressure error term in the velocity estimate is $\nu^{-1}$ [46], whereas we obtain $\gamma^{-1/2}\nu^{-1/2}$ due to the penalization. This is analogous to the effect of using grad-div stabilization in conforming methods, and as we show below, it can dramatically reduce the velocity error when $\nu$ is small and/or the pressure is large. Further, due to (32b) in the limit $\gamma \to \infty$ also the discrete pressure is pressure-robust in the sense that its $L^2$ difference from the discrete approximation of the true solution pressure $\pi_0 p$ does not depend on the true solution pressure $p$, see [29, Remark 4.5].

**Proof of Theorem 3.4:** Testing with $(\boldsymbol{v}_h, q_h) = (\boldsymbol{e}_h^{\boldsymbol{u}}, e_h^p) \in \boldsymbol{V}_h^{\text{div}} \times Q_h$ in Corollary 3.3, inserting the error splitting and reordering leads to

$$\begin{aligned}\nu a_h(\boldsymbol{e}_h^{\boldsymbol{u}}, \boldsymbol{e}_h^{\boldsymbol{u}}) + \gamma j_h(\boldsymbol{e}_h^{\boldsymbol{u}}, \boldsymbol{e}_h^{\boldsymbol{u}}) + b_h(\boldsymbol{\eta}^{\boldsymbol{u}}, e_h^p) - b_h(\boldsymbol{e}_h^{\boldsymbol{u}}, e_h^p) + b_h(\boldsymbol{e}_h^{\boldsymbol{u}}, e_h^p)&\\ = \nu a_h(\boldsymbol{\eta}^{\boldsymbol{u}}, \boldsymbol{e}_h^{\boldsymbol{u}}) + \gamma j_h(\boldsymbol{\eta}^{\boldsymbol{u}}, \boldsymbol{e}_h^{\boldsymbol{u}}) + b_h(\boldsymbol{e}_h^{\boldsymbol{u}}, \eta^p).&\end{aligned} \tag{33}$$

On the left-hand side, due to (31), $\nabla \cdot \boldsymbol{\eta}^{\boldsymbol{u}} = 0$ holds weakly, and thus $b_h(\boldsymbol{\eta}^{\boldsymbol{u}}, e_h^p) = 0$ since also $\boldsymbol{\eta}^{\boldsymbol{u}} \in \boldsymbol{H}_0(\text{div};\Omega)$. The other two mixed terms on the left-hand side cancel each other out. On the right-hand side, we observe that $\boldsymbol{\eta}^{\boldsymbol{u}} \in \boldsymbol{H}_0(\text{div};\Omega)$, and thus $[\![\boldsymbol{\eta}^{\boldsymbol{u}}]\!] \cdot \boldsymbol{n}_F = 0$ for all facets $F \in \mathcal{F}_h$, and so $j_h(\boldsymbol{\eta}^{\boldsymbol{u}}, \boldsymbol{e}_h^{\boldsymbol{u}}) = 0$. For the remaining mixed term on the right-hand side, note that $\nabla \cdot \boldsymbol{e}_h^{\boldsymbol{u}}\big|_K \in \mathbb{P}_{k-1}(K)$ and since $\eta^p$ is orthogonal to $\mathbb{P}_{k-1}$, we obtain

$$b_h(\boldsymbol{e}_h^{\boldsymbol{u}}, \eta^p) = -\int_\Omega \eta^p (\nabla_h \cdot \boldsymbol{e}_h^{\boldsymbol{u}}) \,\mathrm{d}\boldsymbol{x} + \sum_{F\in\mathcal{F}_h} \oint_F \{\!\!\{\eta^p\}\!\!\} ([\![\boldsymbol{e}_h^{\boldsymbol{u}}]\!] \cdot \boldsymbol{n}_F) \,\mathrm{d}\boldsymbol{s} = \sum_{F\in\mathcal{F}_h} \oint_F \{\!\!\{\eta^p\}\!\!\} ([\![\boldsymbol{e}_h^{\boldsymbol{u}}]\!] \cdot \boldsymbol{n}_F) \,\mathrm{d}\boldsymbol{s}.$$

In (33), applying discrete coercivity (21) on the left-hand side and boundedness (17) plus Cauchy–Schwarz on the right-hand side now results in

$$\begin{aligned}\nu C_\sigma |\!|\!| \boldsymbol{e}_h^{\boldsymbol{u}} |\!|\!|_e^2 + \gamma |\boldsymbol{e}_h^{\boldsymbol{u}}|_{\text{nj}}^2 &\leqslant \sqrt{\nu} M |\!|\!| \boldsymbol{\eta}^{\boldsymbol{u}} |\!|\!|_{e,\sharp} \sqrt{\nu} |\!|\!| \boldsymbol{e}_h^{\boldsymbol{u}} |\!|\!|_e \\ &+ \left( \sum_{F\in\mathcal{F}_h} \frac{h_F}{\gamma} \left\| \{\!\!\{\eta^p\}\!\!\} \right\|_{L^2(F)}^2 \right)^{1/2} \left( \sum_{F\in\mathcal{F}_h} \frac{\gamma}{h_F} \left\| [\![\boldsymbol{e}_h^{\boldsymbol{u}}]\!] \cdot \boldsymbol{n}_F \right\|_{L^2(F)}^2 \right)^{1/2} = \mathfrak{T}_1 + \mathfrak{T}_2.\end{aligned}$$

We further estimate these right-hand side terms using Young's inequality:

$$\begin{aligned}\mathfrak{T}_1 &\leqslant \frac{\nu M^2}{2C_\sigma} |\!|\!| \boldsymbol{\eta}^{\boldsymbol{u}} |\!|\!|_{e,\sharp}^2 + \frac{\nu C_\sigma}{2} |\!|\!| \boldsymbol{e}_h^{\boldsymbol{u}} |\!|\!|_e^2, \\ \mathfrak{T}_2 &\leqslant \frac{1}{2} \left( \sum_{F\in\mathcal{F}_h} \frac{h_F}{\gamma} \left\| \{\!\!\{\eta^p\}\!\!\} \right\|_{L^2(F)}^2 \right) + \frac{\gamma}{2} |\boldsymbol{e}_h^{\boldsymbol{u}}|_{\text{nj}}^2.\end{aligned}$$

Thus, the terms involving the velocity discretization error $\boldsymbol{e}_h^{\boldsymbol{u}}$ can be absorbed in the left-hand side. For the term involving the average of $\eta^p$, we use the discrete trace inequality, analogously as for the energy estimate, and obtain

$$\sum_{F\in\mathcal{F}_h} \frac{h_F}{\gamma} \left\| \{\!\!\{\eta^p\}\!\!\} \right\|_{L^2(F)}^2 \leqslant C\gamma^{-1} \left\| \eta^p \right\|_{L^2(\Omega)}^2. \tag{35}$$

Combining all above estimates leads to

$$\nu C_\sigma |\!|\!| \boldsymbol{e}_h^{\boldsymbol{u}} |\!|\!|_e^2 + \gamma |\boldsymbol{e}_h^{\boldsymbol{u}}|_{\text{nj}}^2 \leqslant C\nu |\!|\!| \boldsymbol{\eta}^{\boldsymbol{u}} |\!|\!|_{e,\sharp}^2 + C\gamma^{-1} \left\| \eta^p \right\|_{L^2(\Omega)}^2. \tag{36}$$

Reordering, dropping the positive $\boldsymbol{H}_0(\text{div};\Omega)$-stabilization term on the left-hand side, and taking the square root reveals

$$|\!|\!| \boldsymbol{e}_h^{\boldsymbol{u}} |\!|\!|_e \leqslant C \Big[ |\!|\!| \boldsymbol{\eta}^{\boldsymbol{u}} |\!|\!|_{e,\sharp} + \gamma^{-1/2}\nu^{-1/2} \left\| \eta^p \right\|_{L^2(\Omega)} \Big]. \tag{37}$$

Application of the triangle inequality and the fact that $\|\!|\boldsymbol{\eta}^{\boldsymbol{u}}|\!\|_e \leqslant \|\!|\boldsymbol{\eta}^{\boldsymbol{u}}|\!\|_{e,\sharp}$ gives the claim for the velocity error estimate.

For the pressure estimate, we again invoke the discrete inf-sup condition (15) in $\boldsymbol{W}_h/Q_h$ and the error splitting:

$$\beta^* \|e_h^p\|_{L^2(\Omega)} \leqslant \sup_{\boldsymbol{w}_h \in \boldsymbol{W}_h \setminus \{\boldsymbol{0}\}} \frac{b_h(\boldsymbol{w}_h, e_h^p)}{\|\!|\boldsymbol{w}_h|\!\|_e} = \sup_{\boldsymbol{w}_h \in \boldsymbol{W}_h \setminus \{\boldsymbol{0}\}} \left[\frac{b_h(\boldsymbol{w}_h, \eta^p) - b_h(\boldsymbol{w}_h, p - p_h)}{\|\!|\boldsymbol{w}_h|\!\|_e}\right]. \tag{38}$$

It remains to further estimate the numerator in the last term. Arguing similarly as above, $\nabla \cdot \boldsymbol{w}_h\big|_K \in \mathbb{P}_{k-1}(K)$ as $\boldsymbol{W}_h \subset \boldsymbol{V}_h$, and using the fact that $\eta^p$ is orthogonal to $\mathbb{P}_{k-1}$ yields

$$b_h(\boldsymbol{w}_h, \eta^p) = \sum_{F \in \mathcal{F}_h} \oint_F \{\!\{\eta^p\}\!\} ([\![\boldsymbol{w}_h]\!] \cdot \boldsymbol{n}_F) \, \mathrm{d}\boldsymbol{s} = 0, \tag{39}$$

where the last equality makes use of $[\![\boldsymbol{w}_h]\!] \cdot \boldsymbol{n}_F = 0$ for all $F \in \mathcal{F}_h$. Then, Corollary 3.3 and (17) together with $\boldsymbol{w}_h$ being normal-continuous leads to

$$-b_h(\boldsymbol{w}_h, p - p_h) = -b_h(\boldsymbol{w}_h, \pi_0 p - p_h) = \nu a_h(\boldsymbol{u} - \boldsymbol{u}_h, \boldsymbol{w}_h) \leqslant \nu M \|\!|\boldsymbol{u} - \boldsymbol{u}_h|\!\|_{e,\sharp} \|\!|\boldsymbol{w}_h|\!\|_e. \tag{40}$$

This proves (32b). Equivalence of the $\|\!|\cdot|\!\|_e$- and $\|\!|\cdot|\!\|_{e,\sharp}$-norm, the velocity error estimate (32a) and the triangle inequality conclude the proof. ∎

### *3.3 Numerical illustration of the error estimate*

We now present results of a numerical experiment, in order to illustrate Theorem 3.4. Take as the domain the unit square $\Omega = (0,1)^2$, viscosity $\nu = 10^{-4}$, constant interior penalty parameter $\sigma = 4k^2$, and a third order ($k = 3$) DG method with $\mathbb{P}_3^{\mathrm{dc}}/\mathbb{P}_2^{\mathrm{dc}}$ elements on a structured triangular mesh with $h = \frac{1}{32}$. In order to show that the straight-forward addition of the well-known (broken) grad-div stabilization alone is not sufficient for DG methods, the term

$$\gamma_{\mathrm{gd}} \int_\Omega (\nabla_h \cdot \boldsymbol{u}_h)(\nabla_h \cdot \boldsymbol{v}_h) \, \mathrm{d}\boldsymbol{x} \tag{41}$$

is added to the left-hand side of (19)–(20). Therefore, the parameter $\gamma$ now controls the amount of mass flux penalization, whereas $\gamma_{\mathrm{gd}}$ controls the amount of (broken) divergence penalization. We make use of the high-order finite element library `NGSolve` [47].

The problem we consider is a version of the no-flow problem where the exact solution is chosen to be $\boldsymbol{u} = \boldsymbol{0}$ and $p = \sin(2\pi x + 2\pi y)$. Note that $\int_\Omega p \, \mathrm{d}\boldsymbol{x} = 0$ and the corresponding forcing vector is the gradient field $\boldsymbol{f} = \nabla p$. The results of our experiment are shown in Table 1.

As the mass flux stabilization parameter $\gamma$ increases (with fixed $\gamma_{\mathrm{gd}} = 0$), we observe convergence of the velocity to the no-flow solution. However, for fixed $\gamma = 0$ and increasing $\gamma_{\mathrm{gd}}$, we do not observe an improvement in the velocity error, only in the (broken) divergence error. Especially, broken grad-div stabilization cannot reduce the mass flux error $|\boldsymbol{u}_h|_{\mathrm{nj}}$. The discrete pressure is not significantly influenced by either stabilization.

In contrast to Theorem 3.4 which predicts an error reduction with rate $\gamma^{-1/2}$, Table 1 indicates a better (linear) reduction behaving like $\gamma^{-1}$. The purpose of the next section is to consider the limiting behavior of the method as $\gamma \to \infty$, and resolves this issue of scaling with $\gamma$.

Table 1: Errors for the no-flow problem with the mass flux penalized (controlled by $\gamma$) and broken divergence penalized (controlled by $\gamma_{\mathrm{gd}}$) DG method $\mathbb{P}_3^{\mathrm{dc}}/\mathbb{P}_2^{\mathrm{dc}}$, for $\nu = 10^{-4}$ and using a structured triangular mesh with $h = \frac{1}{32}$.

| $\gamma$ | $\gamma_{\mathrm{gd}}$ | $\lVert \boldsymbol{u} - \boldsymbol{u}_h \rVert_0$ | $\lVert \nabla_h(\boldsymbol{u} - \boldsymbol{u}_h) \rVert_0$ | $\lVert p - p_h \rVert_0$ | $\lVert \nabla_h \cdot \boldsymbol{u}_h \rVert_0$ | $\lvert \boldsymbol{u}_h \rvert_{\mathrm{nj}}$ |
|---|---|---|---|---|---|---|
| 0 | 0 | $3.9 \cdot 10^{-5}$ | $1.3 \cdot 10^{-2}$ | $1.3 \cdot 10^{-5}$ | $9.9 \cdot 10^{-3}$ | $4 \cdot 10^{-3}$ |
| 0.1 | 0 | $3.2 \cdot 10^{-6}$ | $1 \cdot 10^{-3}$ | $1.3 \cdot 10^{-5}$ | $4.5 \cdot 10^{-4}$ | $2.7 \cdot 10^{-4}$ |
| 1 | 0 | $3.5 \cdot 10^{-7}$ | $1.1 \cdot 10^{-4}$ | $1.3 \cdot 10^{-5}$ | $4.8 \cdot 10^{-5}$ | $2.9 \cdot 10^{-5}$ |
| 10 | 0 | $3.5 \cdot 10^{-8}$ | $1.1 \cdot 10^{-5}$ | $1.3 \cdot 10^{-5}$ | $4.8 \cdot 10^{-6}$ | $2.9 \cdot 10^{-6}$ |
| $10^2$ | 0 | $3.5 \cdot 10^{-9}$ | $1.1 \cdot 10^{-6}$ | $1.3 \cdot 10^{-5}$ | $4.8 \cdot 10^{-7}$ | $3 \cdot 10^{-7}$ |
| $10^3$ | 0 | $3.5 \cdot 10^{-10}$ | $1.1 \cdot 10^{-7}$ | $1.3 \cdot 10^{-5}$ | $4.8 \cdot 10^{-8}$ | $3 \cdot 10^{-8}$ |
| 0 | 0.1 | $3.5 \cdot 10^{-5}$ | $1.1 \cdot 10^{-2}$ | $1.4 \cdot 10^{-5}$ | $5.1 \cdot 10^{-5}$ | $2.6 \cdot 10^{-3}$ |
| 0 | 1 | $3.5 \cdot 10^{-5}$ | $1.1 \cdot 10^{-2}$ | $1.4 \cdot 10^{-5}$ | $5.2 \cdot 10^{-6}$ | $2.6 \cdot 10^{-3}$ |
| 0 | 10 | $3.5 \cdot 10^{-5}$ | $1.1 \cdot 10^{-2}$ | $1.4 \cdot 10^{-5}$ | $5.2 \cdot 10^{-7}$ | $2.6 \cdot 10^{-3}$ |
| 0 | $10^2$ | $3.5 \cdot 10^{-5}$ | $1.1 \cdot 10^{-2}$ | $1.4 \cdot 10^{-5}$ | $5.2 \cdot 10^{-8}$ | $2.6 \cdot 10^{-3}$ |
| 0 | $10^3$ | $3.5 \cdot 10^{-5}$ | $1.1 \cdot 10^{-2}$ | $1.4 \cdot 10^{-5}$ | $5.2 \cdot 10^{-9}$ | $2.6 \cdot 10^{-3}$ |

### *3.4 Convergence as $\gamma \to \infty$ to the BDM solution*

We now prove a limiting result for the mass flux penalized DG method as $\gamma \to \infty$. In particular, we will prove that the method converges to a weakly divergence-free BDM solution, with rate $O(\gamma^{-1})$. Since BDM optimally approximates Stokes, this result explains the linear convergence with $\gamma^{-1}$ to the true solution, in the numerical test of the previous subsection.

To begin, we precisely define the $\mathrm{BDM}_k$ space $\boldsymbol{W}_h$ by

$$\boldsymbol{W}_h = \left\{ \boldsymbol{v}_h \in \boldsymbol{H}_0(\mathrm{div};\Omega) \colon \boldsymbol{v}_h\big|_K \in \mathbb{P}_k(K) \right\} = \boldsymbol{V}_h \cap \boldsymbol{H}_0(\mathrm{div};\Omega).$$

The corresponding inf-sup stable FE pair is given by $\boldsymbol{W}_h/Q_h$, also denoted $\mathbb{BDM}_k/\mathbb{P}_{k-1}^{\mathrm{dc}}$. For more information on $\boldsymbol{H}(\mathrm{div})$-FEM, the reader is referred to [50, 48]. Note that the pressure spaces for DG and $\boldsymbol{H}(\mathrm{div})$ methods coincide. A sketch of the local degrees of freedom for both DG and $\boldsymbol{H}(\mathrm{div})$ methods in the 2D case with $k = 3$ is shown in Figure 1.

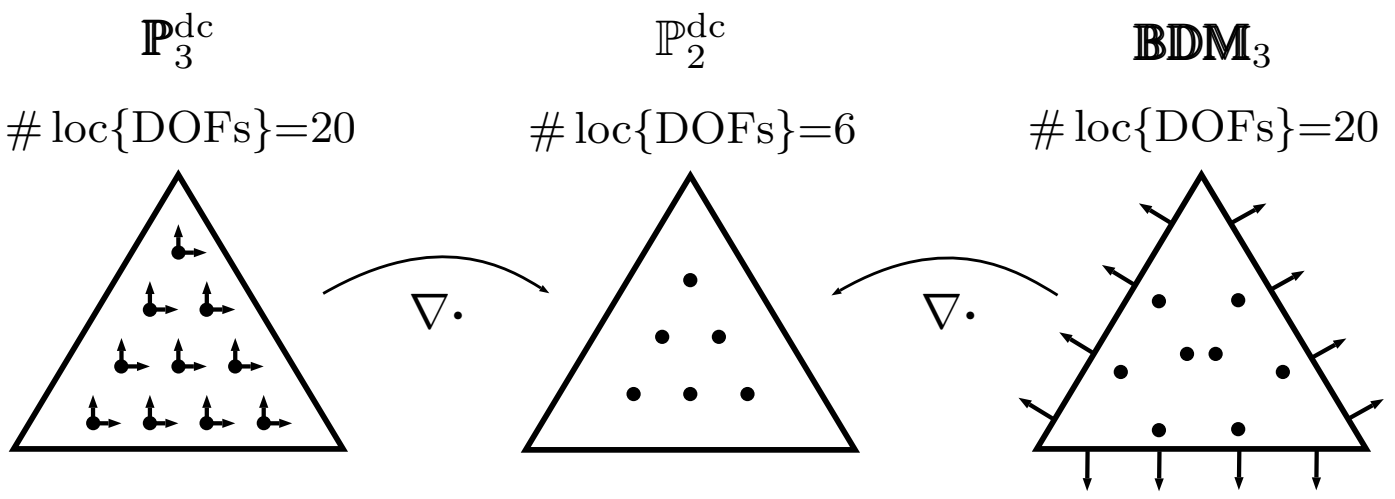

Figure 1: Shown above is a sketch of degrees of freedom in 2D for DG and $\boldsymbol{H}(\mathrm{div})$ methods.

We then introduce the following weakly divergence-free $\boldsymbol{H}(\mathrm{div})$-DG method: Find $(\widehat{\boldsymbol{u}}_h, \widehat{p}_h) \in \boldsymbol{W}_h \times Q_h$ such that for all $(\boldsymbol{w}_h, q_h) \in \boldsymbol{W}_h \times Q_h$,

$$\nu a_h(\widehat{\boldsymbol{u}}_h, \boldsymbol{w}_h) + b(\boldsymbol{w}_h, \widehat{p}_h) - b(\widehat{\boldsymbol{u}}_h, q_h) = (\boldsymbol{f}, \boldsymbol{w}_h). \tag{42}$$

Note that due to the $\boldsymbol{H}(\mathrm{div})$-conformity and in contrast to (19)–(20), the pressure-velocity coupling in (42) is the same as in the continuous weak formulation (2) of the Stokes problem. Since velocities in $\boldsymbol{W}_h$ are normal-continuous, the mass flux penalization naturally vanishes and the SIP bilinear form $a_h$ in (42) acts only on the tangential component of the involved velocities.

The discretely divergence-free subspace of $\boldsymbol{W}_h$ is defined by

$$\boldsymbol{W}_h^{\text{div}} = \{\boldsymbol{w}_h \in \boldsymbol{W}_h : b(\boldsymbol{w}_h, q_h) = 0, \ \forall q_h \in Q_h\},$$

and note that these discretely divergence-free $\boldsymbol{H}(\text{div})$ velocities are even weakly divergence-free, i.e., $\boldsymbol{W}_h^{\text{div}} = \{\boldsymbol{w}_h \in \boldsymbol{H}_0(\text{div};\Omega) : \nabla \cdot \boldsymbol{w}_h = 0\}$.

The aim of this section is to show that the solution $\boldsymbol{u}_h$ of the stabilized DG method (19)–(20) converges to the weakly divergence-free $\boldsymbol{H}(\text{div})$ solution $\widehat{\boldsymbol{u}}_h$ of (42) as $\gamma \to \infty$. To begin, notice that $a_h$ defines a symmetric bilinear form on $\boldsymbol{V}_h$. Due to the inclusion $\boldsymbol{W}_h^{\text{div}} \subset \boldsymbol{V}_h^{\text{div}}$, and since both spaces are finite-dimensional, the orthogonal complement in $\boldsymbol{V}_h^{\text{div}}$, namely

$$\boldsymbol{R}_h^{\perp} = \left(\boldsymbol{W}_h^{\text{div}}\right)^{\perp \boldsymbol{V}_h^{\text{div}}} = \left\{\boldsymbol{v}_h \in \boldsymbol{V}_h^{\text{div}} : a_h(\boldsymbol{v}_h, \boldsymbol{w}_h) = 0, \ \forall \boldsymbol{w}_h \in \boldsymbol{W}_h^{\text{div}}\right\}, \tag{43}$$

makes it possible to obtain the following inner direct sum decomposition:

$$\boldsymbol{V}_h^{\text{div}} = \boldsymbol{W}_h^{\text{div}} \oplus \boldsymbol{R}_h^{\perp}, \quad \boldsymbol{W}_h^{\text{div}} \cap \boldsymbol{R}_h^{\perp} = \{\boldsymbol{0}\}, \quad (\perp \text{ w.r.t. } a_h \text{ inner product}). \tag{44}$$

Thus, $\boldsymbol{W}_h^{\text{div}}$ contains weakly divergence-free, normal-continuous velocities whereas a velocity $\boldsymbol{v}_h^{\perp} \in \boldsymbol{R}_h^{\perp}$ is either $\boldsymbol{0}$ or is not in $\boldsymbol{H}_0(\text{div};\Omega)$, as can be seen as follows. If $\boldsymbol{v}_h^{\perp} \in \boldsymbol{R}_h^{\perp} \subset \boldsymbol{V}_h^{\text{div}}$, then $b_h\left(\boldsymbol{v}_h^{\perp}, q_h\right) = 0$ for all $q_h \in Q_h$. Now, if additionally $\boldsymbol{v}_h^{\perp} \in \boldsymbol{H}_0(\text{div};\Omega)$, then $0 = b_h\left(\boldsymbol{v}_h^{\perp}, q_h\right) = b\left(\boldsymbol{v}_h^{\perp}, q_h\right)$, which implies $\boldsymbol{v}_h \in \boldsymbol{W}_h^{\text{div}}$, $\nabla \cdot \boldsymbol{v}_h^{\perp} = 0$ and thus, due to (44), $\boldsymbol{v}_h^{\perp} = \boldsymbol{0}$. Hence, there are no non-trivial $\boldsymbol{H}_0(\text{div};\Omega)$ velocities in $\boldsymbol{R}_h^{\perp}$.

The following corollary is the key property for showing the convergence of $\boldsymbol{u}_h \to \widehat{\boldsymbol{u}}_h$ as $\gamma \to \infty$.

**Corollary 3.7**

*The mapping* $|\cdot|_{\text{nj}} : \boldsymbol{V}_h \to \mathbb{R}$ *defines a norm on* $\boldsymbol{R}_h^{\perp}$ .

**Proof:** We show that if $\left|\boldsymbol{v}_h^{\perp}\right|_{\text{nj}} = 0$, then $\boldsymbol{v}_h^{\perp} \equiv \boldsymbol{0}$ for all $\boldsymbol{v}_h^{\perp} \in \boldsymbol{R}_h^{\perp}$, which is the only non-trivial property. Let $\boldsymbol{v}_h^{\perp} \in \boldsymbol{R}_h^{\perp}$ with $\left|\boldsymbol{v}_h^{\perp}\right|_{\text{nj}} = 0$. We have already shown above that $\boldsymbol{v}_h^{\perp} \in \boldsymbol{R}_h^{\perp}$ implies that either $\boldsymbol{v}_h^{\perp} \equiv \boldsymbol{0}$ or $\boldsymbol{v}_h^{\perp} \notin \boldsymbol{H}_0(\text{div};\Omega)$. If now $\left|\boldsymbol{v}_h^{\perp}\right|_{\text{nj}} = 0$ and $\boldsymbol{v}_h^{\perp} \notin \boldsymbol{H}_0(\text{div};\Omega)\backslash\{\boldsymbol{0}\}$, then immediately it holds $\boldsymbol{v}_h^{\perp} \equiv \boldsymbol{0}$. This finishes the proof. ■

Next, decompose both the DG solution $\boldsymbol{u}_h = \boldsymbol{u}_h^0 + \boldsymbol{u}_h^{\perp}$ with $\left(\boldsymbol{u}_h^0, \boldsymbol{u}_h^{\perp}\right) \in \boldsymbol{W}_h^{\text{div}} \times \boldsymbol{R}_h^{\perp}$ and the DG test functions $\boldsymbol{v}_h = \boldsymbol{v}_h^0 + \boldsymbol{v}_h^{\perp}$ with $\left(\boldsymbol{v}_h^0, \boldsymbol{v}_h^{\perp}\right) \in \boldsymbol{W}_h^{\text{div}} \times \boldsymbol{R}_h^{\perp}$. Inserting this decomposition into (19) and using the properties of the spaces $\boldsymbol{W}_h^{\text{div}}$ and $\boldsymbol{R}_h^{\perp}$ leads to the following decoupled system:

$$\nu a_h\left(\boldsymbol{u}_h^0, \boldsymbol{v}_h^0\right) = \left(\boldsymbol{f}, \boldsymbol{v}_h^0\right), \quad \forall \boldsymbol{v}_h^0 \in \boldsymbol{W}_h^{\text{div}}, \tag{45a}$$

$$\nu a_h\left(\boldsymbol{u}_h^{\perp}, \boldsymbol{v}_h^{\perp}\right) + \gamma j_h\left(\boldsymbol{u}_h^{\perp}, \boldsymbol{v}_h^{\perp}\right) = \left(\boldsymbol{f}, \boldsymbol{v}_h^{\perp}\right), \quad \forall \boldsymbol{v}_h^{\perp} \in \boldsymbol{R}_h^{\perp}. \tag{45b}$$

Since the solution $\widehat{\boldsymbol{u}}_h$ to the weakly divergence-free $\boldsymbol{H}(\text{div})$-DG is unique, we infer from (45a) that $\boldsymbol{u}_h^0 = \widehat{\boldsymbol{u}}_h$.

Based on this inequality, we can now investigate the limit behavior as $\gamma \to \infty$.

**Theorem 3.8** (*Convergence to divergence-free* $\boldsymbol{H}(\text{div})$ *solution*)

*Let* $\boldsymbol{f} \in \boldsymbol{L}^2(\Omega)$ *and assume that* $\sigma > 0$ *is sufficiently large to guarantee discrete coercivity. Let* $\boldsymbol{u}_h$ *be the solution of the stabilized DG method* (19)–(20) *with* $\gamma > 0$, *and* $\widehat{\boldsymbol{u}}_h$ *be the weakly divergence-free* $\boldsymbol{H}(\text{div})$ *solution of* (42). *Then there exists a constant* $C > 0$, *independent of* $\gamma$, *such that*

$$\left\|\boldsymbol{u}_h - \widehat{\boldsymbol{u}}_h\right\|_e \leqslant C\gamma^{-1} \left\|\boldsymbol{f}\right\|_{\boldsymbol{L}^2(\Omega)}, \tag{46a}$$

$$\left\|p_h - \widehat{p}_h\right\|_{L^2(\Omega)} \leqslant C\gamma^{-1}\nu \left\|\boldsymbol{f}\right\|_{\boldsymbol{L}^2(\Omega)}. \tag{46b}$$

**Remark 3.9:** The generic constant $C$ includes constants arising from norm equivalences of $|\cdot|_{\mathrm{nj}}$ and $|||\cdot|||_e$. While we do not believe this constant depends on $h$, we were unable to prove it. However, we performed several numerical experiments (for brevity not shown here) that did not show any increase in $C$ under mesh refinement.

**Proof of Theorem 3.8:** First note that

$$|||\boldsymbol{u}_h - \widehat{\boldsymbol{u}}_h|||_e = \left|\left|\left|\boldsymbol{u}_h - \boldsymbol{u}_h^0\right|\right|\right|_e = \left|\left|\left|\boldsymbol{u}_h^\perp\right|\right|\right|_e, \tag{47}$$

where $\boldsymbol{u}_h = \boldsymbol{u}_h^0 + \boldsymbol{u}_h^\perp$, $\boldsymbol{u}_h^0 \in \boldsymbol{W}_h^{\mathrm{div}}$, $\boldsymbol{u}_h^\perp \in \boldsymbol{R}_h^\perp$, is the orthogonal decomposition according to (44).

Testing symmetrically with $\boldsymbol{v}_h^\perp = \boldsymbol{u}_h^\perp$ in (45b), using discrete coercivity (21) on the left-hand side and Cauchy–Schwarz on the right-hand side leads to

$$\nu C_\sigma \left|\left|\left|\boldsymbol{u}_h^\perp\right|\right|\right|_e^2 + \gamma\left|\boldsymbol{u}_h^\perp\right|_{\mathrm{nj}}^2 \leqslant \|\boldsymbol{f}\|_{\boldsymbol{L}^2(\Omega)} \left\|\boldsymbol{u}_h^\perp\right\|_{\boldsymbol{L}^2(\Omega)}. \tag{48}$$

Due to the fact that in finite-dimensional spaces all norms are equivalent, an application of Corollary 3.7 on the left-hand side of (48) (simply drop the viscous energy norm multiplied by $\nu$) and Poincaré–Friedrichs (cf., for example, [16, Corollary 5.4]) on the right-hand side leads to

$$\begin{aligned}\gamma C\left|\left|\left|\boldsymbol{u}_h^\perp\right|\right|\right|_e^2 &\leqslant \nu C_\sigma \left|\left|\left|\boldsymbol{u}_h^\perp\right|\right|\right|_e^2 + \gamma\left|\boldsymbol{u}_h^\perp\right|_{\mathrm{nj}}^2 \\ &\leqslant \|\boldsymbol{f}\|_{\boldsymbol{L}^2(\Omega)} \left\|\boldsymbol{u}_h^\perp\right\|_{\boldsymbol{L}^2(\Omega)} \leqslant C_{\mathrm{PF}} \|\boldsymbol{f}\|_{\boldsymbol{L}^2(\Omega)} \left|\left|\left|\boldsymbol{u}_h^\perp\right|\right|\right|_e.\end{aligned}$$

Reordering shows the first claim.

For the pressure convergence, we use the discrete inf-sup condition (15) for the FE pair $\boldsymbol{W}_h/Q_h$. Since both methods use the same pressure space, we can consider $q_h = p_h - \widehat{p}_h \in Q_h$:

$$\begin{aligned}\beta^* \|p_h - \widehat{p}_h\|_{L^2(\Omega)} &\leqslant \sup_{\boldsymbol{w}_h \in \boldsymbol{W}_h \setminus \{\boldsymbol{0}\}} \frac{b_h(\boldsymbol{w}_h, p_h - \widehat{p}_h)}{|||\boldsymbol{w}_h|||_e} \\ &= \sup_{\boldsymbol{w}_h \in \boldsymbol{W}_h \setminus \{\boldsymbol{0}\}} \frac{\nu a_h(\widehat{\boldsymbol{u}}_h - \boldsymbol{u}_h, \boldsymbol{w}_h)}{|||\boldsymbol{w}_h|||_e} \leqslant \nu M |||\widehat{\boldsymbol{u}}_h - \boldsymbol{u}_h|||_{e,\sharp}.\end{aligned}$$

Here, we relied on the properties of $\boldsymbol{H}(\mathrm{div})$ methods and used the boundedness of $a_h$ (17). The final step is to acknowledge that $|||\cdot|||_e$- and $|||\cdot|||_{e,\sharp}$-norms are uniformly equivalent on $\boldsymbol{V}_h$. Reordering and (46a) shows the convergence rate for the pressure. ■

**Remark 3.10:** To avoid the issue of the constant $C$ dependencies in Theorem 3.8, different approaches could be considered that may potentially show the convergence result as $\gamma \to \infty$ without using norm equivalences. One potential strategy could be to use hybridization and the idea that both the $\boldsymbol{L}^2$-DG solution $\boldsymbol{u}_h$ and the $\boldsymbol{H}(\mathrm{div})$-DG solution $\widehat{\boldsymbol{u}}_h$ can equivalently be characterized as solutions of the following problem: Find $(\boldsymbol{u}_h', p_h') \in \boldsymbol{V}_h^{\mathrm{div}} \times Q_h'$ such that for all $(\boldsymbol{v}_h', q_h') \in \boldsymbol{V}_h^{\mathrm{div}} \times Q_h'$,

$$\begin{cases}\nu a_h(\boldsymbol{u}_h', \boldsymbol{v}_h') + d_h(\boldsymbol{v}_h', p_h') = (\boldsymbol{f}, \boldsymbol{v}_h'), \\ d_h(\boldsymbol{u}_h', q_h') - \dfrac{1}{\gamma} e_h(p_h', q_h') = 0.\end{cases} \tag{49}$$

Here, $Q_h'$ is the Lagrange multiplier space (hybrid pressure space) which, loosely speaking, consists of piecewise scalar-valued polynomials on the skeleton $\mathcal{F}_h$. The additional bilinear forms are defined by $d_h(\boldsymbol{v}_h', q_h') = \sum_{F\in\mathcal{F}_h} \oint_F q_h'([\![\boldsymbol{v}_h']\!] \cdot \boldsymbol{n}_F)\,\mathrm{d}\boldsymbol{s}$ and $e_h(p_h', q_h') = \sum_{F\in\mathcal{F}_h} \oint_F p_h' q_h'\,\mathrm{d}\boldsymbol{s}$. Thus, the case $\gamma = \infty$ results in the $\boldsymbol{H}(\mathrm{div})$ solution $\boldsymbol{u}_h' = \widehat{\boldsymbol{u}}_h$ and for any finite $\gamma > 0$, the stabilized $\boldsymbol{L}^2$ solution $\boldsymbol{u}_h' = \boldsymbol{u}_h$ is recovered. Additionally, one obtains $p_h' = \gamma([\![\boldsymbol{u}_h']\!] \cdot \boldsymbol{n}_F)$. The desired convergence as $\gamma \to \infty$ can now be obtained as an application of perturbed saddle point problems; cf., for example, [5, Ch. III, §4, Cor. 4.15]. Verifying the desired assumptions remains to be done, and is not trivial.

We now illustrate Theorem 3.8 with a numerical experiment. We take $\Omega = (0,1)^2$, $\nu = 10^{-3}$, $\sigma = 4k^2$ and the exact solution as

$$\boldsymbol{u} = \begin{pmatrix} \pi \sin^2(\pi x)\sin(2\pi y) \\ -\pi\sin(2\pi x)\sin^2(\pi y) \end{pmatrix}, \quad p = \cos(\pi x)\sin(\pi y).$$

The corresponding right-hand side is

$$\boldsymbol{f} = \begin{pmatrix} -\nu 2\pi^3(2\cos(2\pi x) - 1)\sin(2\pi y) - \pi\sin(\pi x)\sin(\pi y) \\ \nu 2\pi^3\sin(2\pi x)(2\cos(2\pi y) - 1) + \pi\cos(\pi x)\cos(\pi y) \end{pmatrix}.$$

We use a structured triangular mesh with $h = \frac{1}{20}$, $k = 3$, and $(\boldsymbol{u}_h, p_h) \in \boldsymbol{V}_h \times Q_h$. Results (obtained with `NGSolve`) are shown in Table 2. We observe the $O(\gamma^{-1})$ convergence in both pressure and velocity, as $\gamma \to \infty$.

Table 2: Convergence behavior of the $\boldsymbol{H}(\mathrm{div})$-stabilized $\mathbb{P}_3^{\mathrm{dc}}/\mathbb{P}_2^{\mathrm{dc}}$ DG solution $(\boldsymbol{u}_h, p_h) \in \boldsymbol{V}_h \times Q_h$ to the weakly divergence-free $\mathbb{BDM}_3/\mathbb{P}_2^{\mathrm{dc}}$ $\boldsymbol{H}(\mathrm{div})$-DG solution $(\widehat{\boldsymbol{u}}_h, \widehat{p}_h) \in \boldsymbol{W}_h \times Q_h$ as $\gamma \to \infty$.

| $\gamma$ | $\Vert\boldsymbol{u}_h - \widehat{\boldsymbol{u}}_h\Vert_{\boldsymbol{L}^2}$ | $\Vert\nabla_h(\boldsymbol{u}_h - \widehat{\boldsymbol{u}}_h)\Vert_{\boldsymbol{L}^2}$ | $\Vert p_h - \widehat{p}_h\Vert_{L^2}$ |
|---|---|---|---|
| 0 | $1.63 \cdot 10^{-5}$ | $3.01 \cdot 10^{-3}$ | $4.73 \cdot 10^{-6}$ |
| 1 | $1.05 \cdot 10^{-6}$ | $2.05 \cdot 10^{-4}$ | $3.97 \cdot 10^{-7}$ |
| 10 | $1.14 \cdot 10^{-7}$ | $2.23 \cdot 10^{-5}$ | $4.31 \cdot 10^{-8}$ |
| $10^2$ | $1.17 \cdot 10^{-8}$ | $2.28 \cdot 10^{-6}$ | $4.2 \cdot 10^{-9}$ |
| $10^3$ | $1.53 \cdot 10^{-9}$ | $2.8 \cdot 10^{-7}$ | $5.53 \cdot 10^{-10}$ |

## 4 Stabilization of inf-sup stable DG methods on tensor-product meshes

As mentioned in the previous section, on tensor-product elements (quadrilateral and hexagons) using $\mathbb{Q}_k^{\mathrm{dc}}/\mathbb{Q}_{k-1}^{\mathrm{dc}}$, the situation is slightly more involved since $\nabla_h \cdot \mathbb{Q}_k^{\mathrm{dc}} \not\subset \mathbb{Q}_{k-1}^{\mathrm{dc}}$. In order to demonstrate the difference, we repeat the no-flow test from the previous section with $\mathbb{Q}_k^{\mathrm{dc}}/\mathbb{Q}_{k-1}^{\mathrm{dc}}$ elements on a structured mesh consisting of squares. The results (obtained with `NGSolve`) are shown in Table 3.

Interestingly, neither mass flux penalization alone ($\gamma > 0$, $\gamma_{\mathrm{gd}} = 0$), nor broken grad-div stabilization alone ($\gamma_{\mathrm{gd}} > 0$, $\gamma = 0$), is able to improve the pressure-robustness of the considered method, although each of them independently improve the (broken) divergence error. However, when they are added simultaneously ($\gamma > 0$, $\gamma_{\mathrm{gd}} > 0$), Table 3 clearly indicates that the resulting stabilized DG method yields better velocity error (but as expected does not improve the pressure error). In the following, we briefly sketch how this can be shown by the numerical analysis.

The global discrete spaces in this setting are

$$\boldsymbol{V}_h = \left\{\boldsymbol{v}_h \in \boldsymbol{L}^2(\Omega) \colon \boldsymbol{v}_h\big|_K \in \mathbb{Q}_k(K),\ \forall K \in \mathcal{T}_h\right\}, \tag{50a}$$

$$Q_h = \left\{q_h \in L_0^2(\Omega) \colon q_h\big|_K \in \mathbb{Q}_{k-1}(K),\ \forall K \in \mathcal{T}_h\right\}. \tag{50b}$$

Again, we need an $\boldsymbol{H}(\mathrm{div})$-conforming space $\boldsymbol{W}_h \subset \boldsymbol{H}_0(\mathrm{div};\Omega)$ with $\boldsymbol{W}_h \subset \boldsymbol{V}_h$. There again has to be an $\boldsymbol{H}(\mathrm{div})$ interpolator $\boldsymbol{\pi}_h \colon \boldsymbol{V} \to \boldsymbol{W}_h \subset \boldsymbol{V}_h$ which fulfills the commuting diagram property (31). Finally, the pair $\boldsymbol{W}_h/Q_h$ also has to be inf-sup stable.

**Remark 4.1:** Let us give a brief motivation of why these assumptions are reasonable. In the context of $\mathbb{P}_k^{\mathrm{dc}}/\mathbb{P}_{k-1}^{\mathrm{dc}}$ DG methods on simplicial meshes in Section 3, the corresponding divergence-free $\boldsymbol{H}(\mathrm{div})$ FE pair has been characterized as $\mathbb{BDM}_k/\mathbb{P}_{k-1}^{\mathrm{dc}}$. On tensor-product elements the natural idea would be to take into account the Raviart–Thomas (RT) element $\mathbb{RT}_{[k-1]}$ on quads and hexes; cf. [4, Sections 2.4.1 and 2.5]. For this vector-valued element, it is known that $\nabla \cdot \mathbb{RT}_{[k-1]} = \mathbb{Q}_{k-1}^{\mathrm{dc}}$ and, furthermore, $\mathbb{Q}_{k-1}^{\mathrm{dc}} \subset \mathbb{RT}_{[k-1]} \subset \mathbb{Q}_k^{\mathrm{dc}}$. Thus, the validity of the three conditions is ensured.

Table 3: Errors for the no-flow problem with the mass flux stabilization (controlled by $\gamma$) and broken grad-div (controlled by $\gamma_{\mathrm{gd}}$) stabilized DG method $\mathbb{Q}_3^{\mathrm{dc}}/\mathbb{Q}_2^{\mathrm{dc}}$ with $\nu = 10^{-4}$ on a structured quadratic mesh with $h = \frac{1}{32}$.

| $\gamma$ | $\gamma_{\mathrm{gd}}$ | $\|\boldsymbol{u} - \boldsymbol{u}_h\|_0$ | $\|\nabla_h(\boldsymbol{u} - \boldsymbol{u}_h)\|_0$ | $\|p - p_h\|_0$ | $\|\nabla_h \cdot \boldsymbol{u}_h\|_0$ | $\|\boldsymbol{u}_h\|_{\mathrm{nj}}$ |
|---|---|---|---|---|---|---|
| 0 | 0 | $4.3 \cdot 10^{-5}$ | $1.9 \cdot 10^{-2}$ | $3.1 \cdot 10^{-7}$ | $1.8 \cdot 10^{-2}$ | $6.1 \cdot 10^{-3}$ |
| 1 | 0 | $1.6 \cdot 10^{-5}$ | $4.6 \cdot 10^{-3}$ | $3 \cdot 10^{-7}$ | $4.9 \cdot 10^{-5}$ | $1.7 \cdot 10^{-5}$ |
| 10 | 0 | $1.6 \cdot 10^{-5}$ | $4.6 \cdot 10^{-3}$ | $3 \cdot 10^{-7}$ | $4.9 \cdot 10^{-6}$ | $1.7 \cdot 10^{-6}$ |
| $10^2$ | 0 | $1.6 \cdot 10^{-5}$ | $4.6 \cdot 10^{-3}$ | $3 \cdot 10^{-7}$ | $4.9 \cdot 10^{-7}$ | $1.7 \cdot 10^{-7}$ |
| $10^3$ | 0 | $1.6 \cdot 10^{-5}$ | $4.6 \cdot 10^{-3}$ | $3 \cdot 10^{-7}$ | $4.9 \cdot 10^{-8}$ | $1.7 \cdot 10^{-8}$ |
| 0 | 1 | $1.6 \cdot 10^{-5}$ | $5.2 \cdot 10^{-3}$ | $5.5 \cdot 10^{-6}$ | $5.6 \cdot 10^{-6}$ | $1.1 \cdot 10^{-3}$ |
| 0 | 10 | $1.6 \cdot 10^{-5}$ | $5.2 \cdot 10^{-3}$ | $5.5 \cdot 10^{-6}$ | $5.6 \cdot 10^{-7}$ | $1.1 \cdot 10^{-3}$ |
| 0 | $10^2$ | $1.6 \cdot 10^{-5}$ | $5.2 \cdot 10^{-3}$ | $5.5 \cdot 10^{-6}$ | $5.6 \cdot 10^{-8}$ | $1.1 \cdot 10^{-3}$ |
| 0 | $10^3$ | $1.6 \cdot 10^{-5}$ | $5.2 \cdot 10^{-3}$ | $5.5 \cdot 10^{-6}$ | $5.6 \cdot 10^{-9}$ | $1.1 \cdot 10^{-3}$ |
| 1 | 1 | $5.8 \cdot 10^{-6}$ | $1.7 \cdot 10^{-3}$ | $4.8 \cdot 10^{-6}$ | $4.9 \cdot 10^{-6}$ | $1.4 \cdot 10^{-5}$ |
| 10 | 10 | $1.2 \cdot 10^{-6}$ | $3.5 \cdot 10^{-4}$ | $4.8 \cdot 10^{-6}$ | $4.9 \cdot 10^{-7}$ | $2.2 \cdot 10^{-6}$ |
| $10^2$ | $10^2$ | $1.4 \cdot 10^{-7}$ | $4 \cdot 10^{-5}$ | $4.8 \cdot 10^{-6}$ | $4.9 \cdot 10^{-8}$ | $2.4 \cdot 10^{-7}$ |
| $10^3$ | $10^3$ | $1.4 \cdot 10^{-8}$ | $4.1 \cdot 10^{-6}$ | $4.8 \cdot 10^{-6}$ | $4.9 \cdot 10^{-9}$ | $2.5 \cdot 10^{-8}$ |

For the sake of brevity in the analysis, we only allow one stabilization parameter $\gamma$ for both mass flux and broken grad-div stabilization and redefine the stabilization bilinear form as

$$j_h(\boldsymbol{w}_h, \boldsymbol{v}_h) = \int_\Omega (\nabla_h \cdot \boldsymbol{w}_h)(\nabla_h \cdot \boldsymbol{v}_h)\,\mathrm{d}\boldsymbol{x} + \sum_{F \in \mathcal{F}_h} \frac{1}{h_F} \oint_F ([\![\boldsymbol{w}_h]\!] \cdot \boldsymbol{n}_F)([\![\boldsymbol{v}_h]\!] \cdot \boldsymbol{n}_F)\,\mathrm{d}\boldsymbol{s}. \tag{51}$$

**Theorem 4.2 (*Error estimate on tensor-product meshes*)**

*Under the assumptions of Theorem 3.4, the following error estimate holds true for the stabilized DG method on tensor-product elements:*

$$\|\!|\boldsymbol{u} - \boldsymbol{u}_h|\!\|_e \leqslant C\Big[\|\!|\boldsymbol{\eta}^{\boldsymbol{u}}|\!\|_{e,\sharp} + \gamma^{-1/2}\nu^{-1/2}\,\|\eta^p\|_{L^2(\Omega)}\Big], \tag{52a}$$

$$\|p - p_h\| \leqslant C\nu\|\!|\boldsymbol{\eta}^{\boldsymbol{u}}|\!\|_{e,\sharp} + \left[C\sqrt{\frac{\nu}{\gamma}} + 2\right]\|\eta^p\|_{L^2(\Omega)}\,. \tag{52b}$$

**Proof:** The proof is very similar to that of Theorem 3.4. We only comment on the parts where differences occur between simplicial and tensor-product elements.

First, we observe that also for the new stabilization term, it holds $j_h(\boldsymbol{\eta}^{\boldsymbol{u}}, \boldsymbol{e}_h^{\boldsymbol{u}}) = 0$ since $\boldsymbol{\eta}^{\boldsymbol{u}} \in \boldsymbol{H}_0(\mathrm{div};\Omega)$ and $\nabla \cdot \boldsymbol{\eta}^{\boldsymbol{u}} = 0$ hold exactly. The main difference occurs in the treatment of the mixed term of $b_h$, where, after applying Cauchy–Schwarz, Young, and trace inequalities, we obtain

$$\begin{aligned} b_h(\boldsymbol{e}_h^{\boldsymbol{u}}, \eta^p) &= -\int_\Omega \eta^p(\nabla_h \cdot \boldsymbol{e}_h^{\boldsymbol{u}})\,\mathrm{d}\boldsymbol{x} + \sum_{F \in \mathcal{F}_h} \oint_F \{\!\{\eta^p\}\!\}([\![\boldsymbol{e}_h^{\boldsymbol{u}}]\!] \cdot \boldsymbol{n}_F)\,\mathrm{d}\boldsymbol{s} \\ &\leqslant C\gamma^{-1}\,\|\eta^p\|_{L^2(\Omega)}^2 + \frac{\gamma}{2}\,\|\nabla_h \cdot \boldsymbol{e}_h^{\boldsymbol{u}}\|_{L^2(\Omega)}^2 + \frac{\gamma}{2}|\boldsymbol{e}_h^{\boldsymbol{u}}|_{\mathrm{nj}}^2. \end{aligned}$$

Applying the same ideas as in the proof for simplices, we arrive at

$$\nu C_\sigma \|\boldsymbol{e}_h^{\boldsymbol{u}}\|_e^2 + \gamma |\boldsymbol{e}_h^{\boldsymbol{u}}|_{\mathrm{nj}}^2 + \gamma \left\| \nabla_h \cdot \boldsymbol{e}_h^{\boldsymbol{u}} \right\|_{L^2(\Omega)}^2 \leqslant C\nu \|\boldsymbol{\eta}^{\boldsymbol{u}}\|_{e,\sharp}^2 + C\gamma^{-1} \left\| \eta^p \right\|_{L^2(\Omega)}^2, \tag{53}$$

which yields the velocity error estimate. It is critical to note here that without the (broken) grad-div stabilization, the term $C\gamma^{-1} \|\eta^p\|_{L^2(\Omega)}^2$ on the right-hand side of (53) would be $C\nu^{-1} \|\eta^p\|_{L^2(\Omega)}^2$.

For the pressure estimate, the discrete inf-sup condition (15) in $\boldsymbol{W}_h/Q_h$ is again essential. However, since $\nabla_h \cdot \mathbb{Q}_k^{\mathrm{dc}} \not\subset \mathbb{Q}_{k-1}^{\mathrm{dc}}$, for $\boldsymbol{w}_h \in \boldsymbol{W}_h$ we have to estimate

$$b_h(\boldsymbol{w}_h, \eta^p) = -\int_\Omega \eta^p (\nabla \cdot \boldsymbol{w}_h) \,\mathrm{d}\boldsymbol{x} \leqslant \|\eta^p\|_{L^2(\Omega)} \|\nabla \cdot \boldsymbol{w}_h\|_{L^2(\Omega)}. \tag{54}$$

Now, $\|\nabla \cdot \boldsymbol{w}_h\|_{L^2(\Omega)} \leqslant \|\boldsymbol{w}_h\|_e$ leads immediately to the same situation as in the proof of Theorem 3.4. ■

**Remark 4.3:** Concerning the convergence of the stabilized $\boldsymbol{V}_h/Q_h$ DG solution to an $\boldsymbol{H}(\mathrm{div})$ solution as $\gamma \to \infty$, basically the same arguments as in Section 3.2 can be applied. Unfortunately, while considering the $\boldsymbol{H}(\mathrm{div})$ FE pair $\boldsymbol{W}_h/Q_h = \mathbb{RT}_{[k-1]}/\mathbb{Q}_{k-1}^{\mathrm{dc}}$ is sufficient for the error analysis leading to Theorem 4.2, it has the shortcoming that $\boldsymbol{W}_h \neq \boldsymbol{V}_h \cap \boldsymbol{H}_0(\mathrm{div};\Omega)$ (in fact, the major problem is that $\boldsymbol{W}_h$ has only $(k-1)$th-order approximation properties). For a convergence result as $\gamma \to \infty$, it is essential to have a space $\boldsymbol{W}_h^+ \subset \boldsymbol{H}_0(\mathrm{div};\Omega)$ fulfilling $\boldsymbol{W}_h^+ = \boldsymbol{V}_h \cap \boldsymbol{H}_0(\mathrm{div};\Omega)$. Such a space $\boldsymbol{W}_h^+$ could be constructed, but, to the best of the authors' knowledge, has not been previously used in the literature. Thus, we do not go into detail at this point. However, let us mention that one can show that $\mathbb{RT}_{[k-1]} \subset \boldsymbol{W}_h^+ \subset \mathbb{Q}_k^{\mathrm{dc}}$.

**Remark 4.4:** In light of the problem of explicitly having the space $\boldsymbol{W}_h^+$ available, we note that it is nonetheless possible to obtain a discrete solution $(\widehat{\boldsymbol{u}}_h, \widehat{p}_h) \in \boldsymbol{W}_h^+ \times Q_h$. To this end, a Lagrange multiplier technique can be used where additional unknowns on facets are introduced which represent the normal component of the velocity. Similar ideas are frequently employed in hybridized discontinuous Galerkin (HDG) methods. For more information, we refer to [32, Remark 7] and/or [10, 45].

## 5 Improving pressure-robustness in Crouzeix–Raviart approximations

The mass flux penalization can also be applied, with similar results, to the Crouzeix–Raviart (CR) element, and we include some results for completeness. In a sense, it is somewhat easier to analyze than the DG case considered above, and we find a similar fundamental result: penalization of the mass flux reduces the effect of the pressure error on the velocity error. Such a result is proven implicitly for CR elements in a recent work of Burman and Hansbo [7] for the Darcy–Stokes problem, where multiple stabilizations were used. Interestingly, the motivation for using the stabilization in that work was 'to control the nonconformity emanating from the pressure term', and in effect they proved something similar to what we prove above for DG: the scaling of the pressure term in the error estimate is improved from $\nu^{-1}$ to $\nu^{-1/2}\gamma^{-1/2}$.

However, there is (seemingly) a potential negative consequence for CR that is not an issue with DG: the use of the stabilization seemingly increases the scaling of the velocity error in the error estimate. Hence it is unclear whether large stabilization parameters can be used without negative consequence, as can be done in the DG case (at least, up to difficulties in linear solvers due to matrix conditioning).

The nonconforming Crouzeix–Raviart finite element velocity and pressure spaces are defined by:

$$\begin{aligned} \boldsymbol{V}_h &= \left\{ \boldsymbol{v}_h \in \boldsymbol{L}^2(\Omega) \colon \boldsymbol{v}_h\big|_K \in \mathbb{P}_1(K),\ \forall K \in \mathcal{T}_h,\ \oint_F [\![\boldsymbol{v}_h]\!] = \boldsymbol{0},\ \forall F \in \mathcal{F}_h \right\}, \\ Q_h &= \left\{ q_h \in L_0^2(\Omega) \colon q_h\big|_K \in \mathbb{P}_0(K),\ \forall K \in \mathcal{T}_h \right\}. \end{aligned}$$

The discrete bilinear forms for CR elements are defined as

$$a_h(\boldsymbol{w}_h, \boldsymbol{v}_h) = \int_\Omega \nabla_h \boldsymbol{w}_h : \nabla_h \boldsymbol{v}_h \, \mathrm{d}\boldsymbol{x},$$
$$b_h(\boldsymbol{w}_h, q_h) = -\int_\Omega q_h (\nabla_h \cdot \boldsymbol{w}_h) \, \mathrm{d}\boldsymbol{x}.$$

$\boldsymbol{V}_h$ is equipped with the norm

$$\|\boldsymbol{v}_h\|_{1,h} := \left( \int_\Omega \nabla_h \boldsymbol{v}_h : \nabla_h \boldsymbol{v}_h \, \mathrm{d}\boldsymbol{x} \right)^{1/2} = \left( \sum_{K \in \mathcal{T}_h} \|\nabla \boldsymbol{v}_h\|^2_{\boldsymbol{L}^2(K)} \right)^{1/2}.$$

We propose to consider CR together with the mass flux penalization: find $(\boldsymbol{u}_h, \boldsymbol{v}_h) \in \boldsymbol{V}_h \times Q_h$ such that for all $(\boldsymbol{v}_h, q_h) \in \boldsymbol{V}_h \times Q_h$

$$\gamma j_h(\boldsymbol{u}_h, \boldsymbol{v}_h) + \nu a_h(\boldsymbol{u}_h, \boldsymbol{v}_h) + b_h(\boldsymbol{v}_h, p_h) = (\boldsymbol{f}, v_h) \tag{55}$$
$$-b_h(\boldsymbol{u}_h, q_h) = 0. \tag{56}$$

Since $j_h(\boldsymbol{u}_h, \boldsymbol{u}_h) \geqslant 0$, the classical well-posedness results for CR (see e.g. [6]) will hold also for (55)–(56) with any fixed $\gamma \geqslant 0$.

We now present an error estimate, which follows from the recent work of Burman and Hansbo for a Darcy–Stokes problem with multiple stabilizations [7], and making the appropriate simplifications (hence we omit the proof).

**Theorem 5.1 (*Error estimate for Crouzeix–Raviart elements*)**

*Let $(\boldsymbol{u}, p) \in \big(\boldsymbol{V} \cap \boldsymbol{H}^2(\Omega)\big) \times \big(Q \cap H^1(\Omega)\big)$ be the Stokes solution, and $(\boldsymbol{u}_h, p_h)$ the solution to (55)–(56) with parameter $\gamma > 0$. Then it holds that*

$$\|\boldsymbol{u} - \boldsymbol{u}_h\|_{1,h} \leqslant Ch\Big[\Big(1 + \gamma^{1/2}\nu^{-1/2}\Big)|\boldsymbol{u}|_{\boldsymbol{H}^2(\Omega)} + \gamma^{-1/2}\nu^{-1/2}|p|_{H^1(\Omega)}\Big].$$

**Remark 5.2:** This error estimates reveal, just as in the DG case, that the facet jump stabilization reduces the negative effect of the pressure on the velocity error by changing the coefficient of the pressure from $\nu^{-1}$ to $\nu^{-1/2}\gamma^{-1/2}$. However, the dependence on $\nu^{-1/2}$ in the term $\big(1 + \gamma^{1/2}\nu^{-1/2}\big)|\boldsymbol{u}|_{\boldsymbol{H}^2(\Omega)}$ seems to be unavoidable for the nonconforming CR-element, indicating the danger of over-stabilization, which is nothing more than a kind of Poisson locking for the divergence-free limit as $\gamma \to \infty$.

## 6 Application to more complex flows

We show here that the pressure-robust approach above can have a considerable impact on more complicated problems than the steady incompressible Stokes equations. We consider first a numerical test for steady incompressible Navier–Stokes equations (NSE), followed by a numerical test for non-isothermal flow.

### *6.1 Steady Navier–Stokes equations*

Consider now the steady incompressible Navier–Stokes equations

$$\begin{cases} -\nu\Delta\boldsymbol{u} + (\boldsymbol{u} \cdot \nabla)\boldsymbol{u} + \nabla p = \boldsymbol{f} & \text{in } \Omega, \\ \nabla \cdot \boldsymbol{u} = 0 & \text{in } \Omega, \end{cases} \tag{57}$$

with inhomogeneous Dirichlet velocity boundary condition $\boldsymbol{g}_{\mathrm{D}}$. The new approximately pressure-robust space discretization approach will be superior, whenever a difficult pressure spoils the velocity error.

In order to show this, we construct a simple polynomial potential flow $\boldsymbol{u} = \nabla\chi$ on the unit square $\Omega = (-1,1)^2$, where $\chi$ is defined by the real part of the complex polynomial $z^5$, i.e., $\chi = \mathrm{Re}(z^5) = \mathrm{Re}(x+iy)^5 = x^5 - 10x^3y^2 + 5xy^4$. The resulting flow consists of ten colliding jets, meeting at the stagnation point $(0,0)$, see Figure 2. Indeed, $(\boldsymbol{u}, p) = (\nabla\chi, -\frac{1}{2}|\nabla\chi|^2)$ solves the steady incompressible NSE (57) for $\boldsymbol{f} = \boldsymbol{0}$ with appropriate inhomogeneous Dirichlet velocity boundary conditions, for all $\nu > 0$ [37]. Since $p$ is a polynomial of order eight, a pressure-robust low-order DG space discretization should be comparably accurate — at least at non-negligible Reynolds numbers — as a high-order DG approximation which is not pressure-robust.

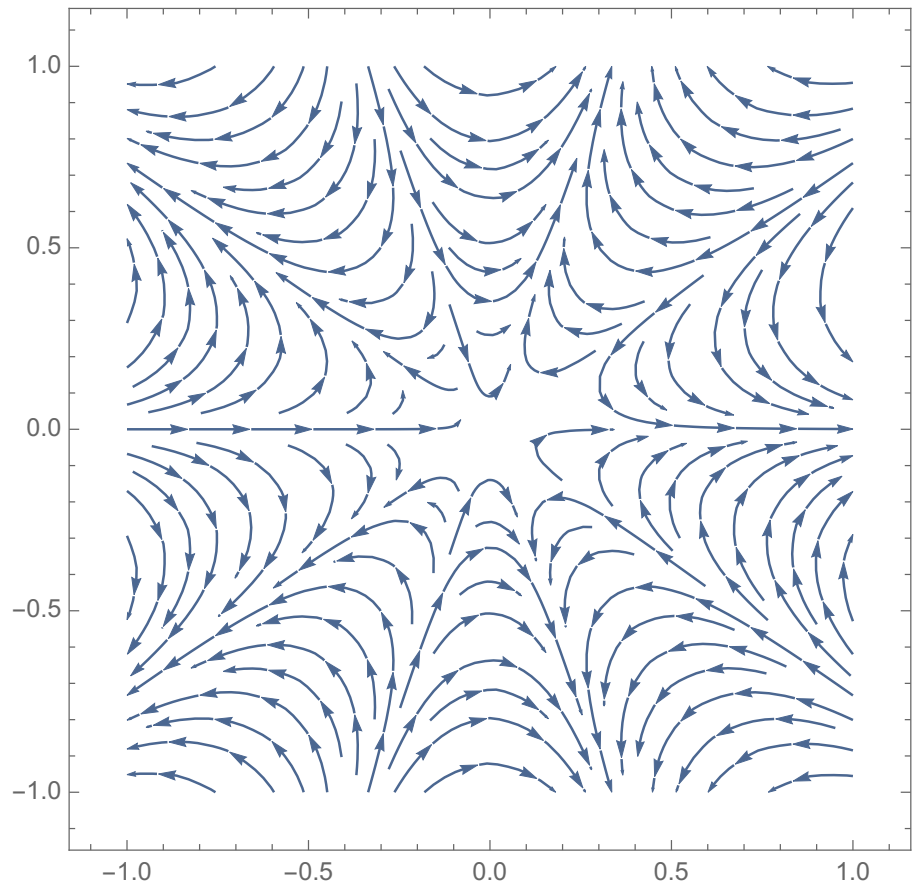

Figure 2: Streamlines for the steady incompressible Navier–Stokes benchmark with ten jets and a stagnation point at $(0,0)$.

The proposed stabilized DG method for the stationary Navier–Stokes equations (57) is the following: Find $(\boldsymbol{u}_h, p_h) \in \boldsymbol{V}_h \times Q_h$ satisfying $\forall\,(\boldsymbol{v}_h, q_h) \in \boldsymbol{V}_h \times Q_h$,

$$\begin{aligned}\nu a_h(\boldsymbol{u}_h, \boldsymbol{v}_h) + c_h(\boldsymbol{u}_h; \boldsymbol{u}_h, \boldsymbol{v}_h) + b_h(\boldsymbol{v}_h, p_h) + \gamma j_h(\boldsymbol{u}_h, \boldsymbol{v}_h) &= \nu a_h^{\partial}(\boldsymbol{g}_{\mathrm{D}}; \boldsymbol{v}_h) + \gamma j_h^{\partial}(\boldsymbol{g}_{\mathrm{D}}; \boldsymbol{v}_h),\\ b_h(\boldsymbol{u}_h, q_h) &= b_h^{\partial}(\boldsymbol{g}_{\mathrm{D}}; q_h).\end{aligned}$$

Due to the weak imposition of a non-zero Dirichlet boundary condition $\boldsymbol{g}_{\mathrm{D}} = \boldsymbol{u}\big|_{\partial\Omega}$, additional boundary facet terms arise from the SIP bilinear form, the pressure-velocity coupling and the mass flux stabilization:

$$\begin{aligned}a_h^{\partial}(\boldsymbol{g}_{\mathrm{D}}; \boldsymbol{v}_h) &= \sum_{F\in\mathcal{F}_h^{\partial}} \frac{\sigma}{h_F} \oint_F \boldsymbol{g}_{\mathrm{D}} \cdot \boldsymbol{v}_h \,\mathrm{d}\boldsymbol{s} - \sum_{F\in\mathcal{F}_h^{\partial}} \oint_F \boldsymbol{g}_{\mathrm{D}} \cdot \nabla\boldsymbol{v}_h\boldsymbol{n} \,\mathrm{d}\boldsymbol{s},\\ b_h^{\partial}(\boldsymbol{g}_{\mathrm{D}}; q_h) &= \sum_{F\in\mathcal{F}_h^{\partial}} \oint_F q_h(\boldsymbol{g}_{\mathrm{D}} \cdot \boldsymbol{n}) \,\mathrm{d}\boldsymbol{s},\\ j_h^{\partial}(\boldsymbol{g}_{\mathrm{D}}; \boldsymbol{v}_h) &= \sum_{F\in\mathcal{F}_h^{\partial}} \frac{1}{h_F} \oint_F (\boldsymbol{g}_{\mathrm{D}} \cdot \boldsymbol{n})(\boldsymbol{v}_h \cdot \boldsymbol{n}) \,\mathrm{d}\boldsymbol{s}.\end{aligned}$$

Again, $\sigma = 4k^2$ is chosen. In contrast to the Stokes problem, we also additionally have to deal with nonlinear inertia effects. For treating this, we choose the following standard convection term with upwinding [16]:

$$c_h(\boldsymbol{w}_h; \boldsymbol{u}_h, \boldsymbol{v}_h) = \int_\Omega (\boldsymbol{w}_h \cdot \nabla_h)\boldsymbol{u}_h \cdot \boldsymbol{v}_h \,\mathrm{d}\boldsymbol{x}$$
$$- \sum_{F \in \mathcal{F}_h^i} \oint_F (\{\!\{\boldsymbol{w}_h\}\!\} \cdot \boldsymbol{n}_F)[\![\boldsymbol{u}_h]\!] \cdot \{\!\{\boldsymbol{v}_h\}\!\} \,\mathrm{d}\boldsymbol{s} + \sum_{F \in \mathcal{F}_h^i} \oint_F \frac{1}{2} |\{\!\{\boldsymbol{w}_h\}\!\} \cdot \boldsymbol{n}_F| [\![\boldsymbol{u}_h]\!] \cdot [\![\boldsymbol{v}_h]\!] \,\mathrm{d}\boldsymbol{s}.$$

For the first numerical test, we use $\nu = 10^{-2}$ and compute using the FE pair $\mathbb{P}_4^{\mathrm{dc}}/\mathbb{P}_3^{\mathrm{dc}}$ on an unstructured triangular mesh with $h = 0.1$. In this case, the penalization term $j_h$ only involves the mass flux contribution, see (18). Following [46], the nonlinear system is linearized using Picard's iteration (with tolerance $10^{-8}$) and for the resulting linear systems, a standard reiteration scheme ensures accurate solutions also for large stabilization parameters $\gamma$. The results can be seen in Table 4, and show a clear and remarkable effect on the velocity error. Since the true velocity solution is a fourth degree polynomial, and the method becomes pressure-robust as $\gamma \to \infty$, we observe very small velocity errors for the larger penalization parameters.

Table 4: Errors for the Navier–Stokes potential flow problem with the mass flux penalized (controlled by $\gamma$) DG method $\mathbb{P}_4^{\mathrm{dc}}/\mathbb{P}_3^{\mathrm{dc}}$ with $\nu = 10^{-2}$ on an unstructured triangular mesh with $h = 0.1$.

| $\gamma$ | $\Vert\boldsymbol{u}-\boldsymbol{u}_h\Vert_0$ | $\Vert\nabla_h(\boldsymbol{u}-\boldsymbol{u}_h)\Vert_0$ | $\Vert p-p_h\Vert_0$ | $\Vert\nabla_h\cdot\boldsymbol{u}_h\Vert_0$ | $\vert\boldsymbol{u}_h\vert_{\mathrm{nj}}$ |
|---|---|---|---|---|---|
| 0 | $1.4\cdot10^{-4}$ | $1.2\cdot10^{-2}$ | $9.1\cdot10^{-4}$ | $9.9\cdot10^{-3}$ | $1.5\cdot10^{-3}$ |
| 1 | $3.2\cdot10^{-5}$ | $3.3\cdot10^{-3}$ | $5.2\cdot10^{-4}$ | $2.6\cdot10^{-3}$ | $6.9\cdot10^{-4}$ |
| 10 | $6.6\cdot10^{-6}$ | $6.6\cdot10^{-4}$ | $4.7\cdot10^{-4}$ | $4.6\cdot10^{-4}$ | $1.3\cdot10^{-4}$ |
| $10^2$ | $7.7\cdot10^{-7}$ | $7.6\cdot10^{-5}$ | $4.7\cdot10^{-4}$ | $5.1\cdot10^{-5}$ | $1.5\cdot10^{-5}$ |
| $10^3$ | $7.6\cdot10^{-8}$ | $7.7\cdot10^{-6}$ | $4.7\cdot10^{-4}$ | $5.2\cdot10^{-6}$ | $1.5\cdot10^{-6}$ |
| $10^4$ | $6.2\cdot10^{-9}$ | $7.6\cdot10^{-7}$ | $4.7\cdot10^{-4}$ | $5.2\cdot10^{-7}$ | $1.5\cdot10^{-7}$ |

We now repeat this test on a tensor-product mesh consisting of quadrilaterals with $h = 0.1$. The corresponding FE pair is thus $\mathbb{Q}_4^{\mathrm{dc}}/\mathbb{Q}_3^{\mathrm{dc}}$ and the proposed stabilized method now also contains the broken grad-div stabilization; i.e., $j_h$ is defined by (51) and the parameter $\gamma$ multiplies both contributions. Table 5 shows the results of this test and indeed, qualitatively the same behavior as on simplicial meshes can be observed.

Table 5: Errors for the Navier–Stokes potential flow problem with the both mass flux penalized and broken grad-div stabilized (both controlled by $\gamma$) DG method $\mathbb{Q}_4^{\mathrm{dc}}/\mathbb{Q}_3^{\mathrm{dc}}$ with $\nu = 10^{-2}$ on a quadrilateral mesh with $h = 0.1$.

| $\gamma$ | $\Vert\boldsymbol{u}-\boldsymbol{u}_h\Vert_0$ | $\Vert\nabla_h(\boldsymbol{u}-\boldsymbol{u}_h)\Vert_0$ | $\Vert p-p_h\Vert_0$ | $\Vert\nabla_h\cdot\boldsymbol{u}_h\Vert_0$ | $\vert\boldsymbol{u}_h\vert_{\mathrm{nj}}$ |
|---|---|---|---|---|---|
| 0 | $1.3\cdot10^{-3}$ | $1.9\cdot10^{-1}$ | $1.1\cdot10^{-2}$ | $1.9\cdot10^{-1}$ | $1.1\cdot10^{-2}$ |
| 10 | $2.4\cdot10^{-6}$ | $2.1\cdot10^{-4}$ | $7.8\cdot10^{-4}$ | $5.7\cdot10^{-5}$ | $2.8\cdot10^{-5}$ |
| $10^2$ | $1.3\cdot10^{-6}$ | $1\cdot10^{-4}$ | $7.8\cdot10^{-4}$ | $5.8\cdot10^{-6}$ | $3\cdot10^{-6}$ |
| $10^3$ | $2.5\cdot10^{-7}$ | $1.8\cdot10^{-5}$ | $7.8\cdot10^{-4}$ | $5.9\cdot10^{-7}$ | $3.2\cdot10^{-7}$ |
| $10^4$ | $3\cdot10^{-8}$ | $2.1\cdot10^{-6}$ | $7.8\cdot10^{-4}$ | $6\cdot10^{-8}$ | $3.2\cdot10^{-8}$ |
| $10^5$ | $5.6\cdot10^{-9}$ | $4.7\cdot10^{-7}$ | $7.8\cdot10^{-4}$ | $6\cdot10^{-9}$ | $3.2\cdot10^{-9}$ |

### *6.2 Non-isothermal flows*

We consider for our final test a differentially heated cavity with infinite Prandtl number (which corresponds physically to silicon oil) and Rayleigh number $Ra = 10^6$. Here we test the CR approximation, with and without mass flux penalization. The problem setup is taken from [28], and we approximate the problem

$$\nabla p - \nu \Delta \boldsymbol{u} = RaT\mathbf{e}_2,$$
$$\nabla \cdot \boldsymbol{u} = 0,$$
$$\boldsymbol{u} \cdot \nabla T - \Delta T = 0,$$

on $\Omega = (0,1)^2$, with $\nu = 1$, no-slip boundary conditions on all walls, insulated boundary conditions on the top and bottom: $\nabla T \cdot \boldsymbol{n}\big|_{\Gamma_T \cup \Gamma_B} = 0$, and Dirichlet temperature conditions on the left and right walls, $T\big|_{\Gamma_L} = 1$, $T\big|_{\Gamma_R} = 0$. The system is approximated using CR elements for the velocity and pressure, and $S_h = \mathbb{P}_1 \cap H^1(\Omega)$ for the temperature approximation. The scheme takes the form

$$\gamma j_h(\boldsymbol{u}_h, \boldsymbol{v}_h) + \nu a_h(\boldsymbol{u}_h, \boldsymbol{v}_h) + b_h(\boldsymbol{v}_h, p_h) = Ra(T_h\mathbf{e}_2, \boldsymbol{v}_h), \tag{58}$$
$$b_h(\boldsymbol{u}_h, q_h) = 0, \tag{59}$$
$$(\boldsymbol{u}_h \cdot \nabla T_h, s_h) + (\nabla T_h, \nabla s_h) = 0, \tag{60}$$

together with appropriate boundary conditions, and Newton's method is used for the arising nonlinearity.

Simulations were run using (58)–(60), for $\gamma \in \{0, 0.1, 10, 10^3\}$ on a 48×48 uniform mesh that was additionally refined once in all cells touching the boundary. Results are shown in Figure 3, and additionally we show a reference solution found using $(\mathbb{P}_2, \mathbb{P}_1, \mathbb{P}_2)$ elements on a 64×64 mesh (we note this solution matches that found in [28]). From the plots, we observe that with no stabilization, the CR solution is very poor: the velocity streamlines are very under-resolved, and the predicted temperature contours are visibly inaccurate. Clear improvement is seen in the solution as the stabilization parameter is increased, and with $\gamma = 10^3$, the CR solution temperature and velocity have essentially converged to that of the reference solution. The pressure is not converged with this largest $\gamma$, although it is significantly improved compared to lower values of $\gamma$ and especially compared to the unstabilized solution.

## 7 Conclusions and future directions

We have analyzed and tested a penalization(s) for nonconforming methods that has an effect on solutions analogous to that of grad-div stabilization for conforming methods. The penalization is a mass flux penalization, and also a broken divergence stabilization in elements where the divergence of the velocity space is not contained in the pressure space.

Our new theoretical results include error estimates that reduce the scaling of the velocity error by the pressure from $\nu^{-1}$ to $\gamma^{-1/2}\nu^{-1/2}$, analogous to the effect of grad-div for conforming methods [28] for the penalized DG method on tensor-product meshes. Additionally, we prove limiting behavior results for the penalized DG method on simplex meshes, and in particular that the limit as $\gamma \to \infty$ is the associated (optimal) BDM solution. Therefore, over-stabilization is not possible in this DG context, surprisingly. We also consider the limiting behavior on tensor-product meshes, and find a limit that is seemingly not in the literature, but is in between successive RT element spaces (and is thus optimal).

For future research directions, since there is now an overwhelming body of work done for grad-div stabilization, we expect that many similar results will also hold with the mass flux penalization for DG methods. This includes further error analyses on more complex problems, linear algebra considerations, turbulence models, and so on. Extension of the limiting behavior results herein to the case of Navier–Stokes equations is also an important next step, as is an extension to more general meshes that allow hanging nodes.

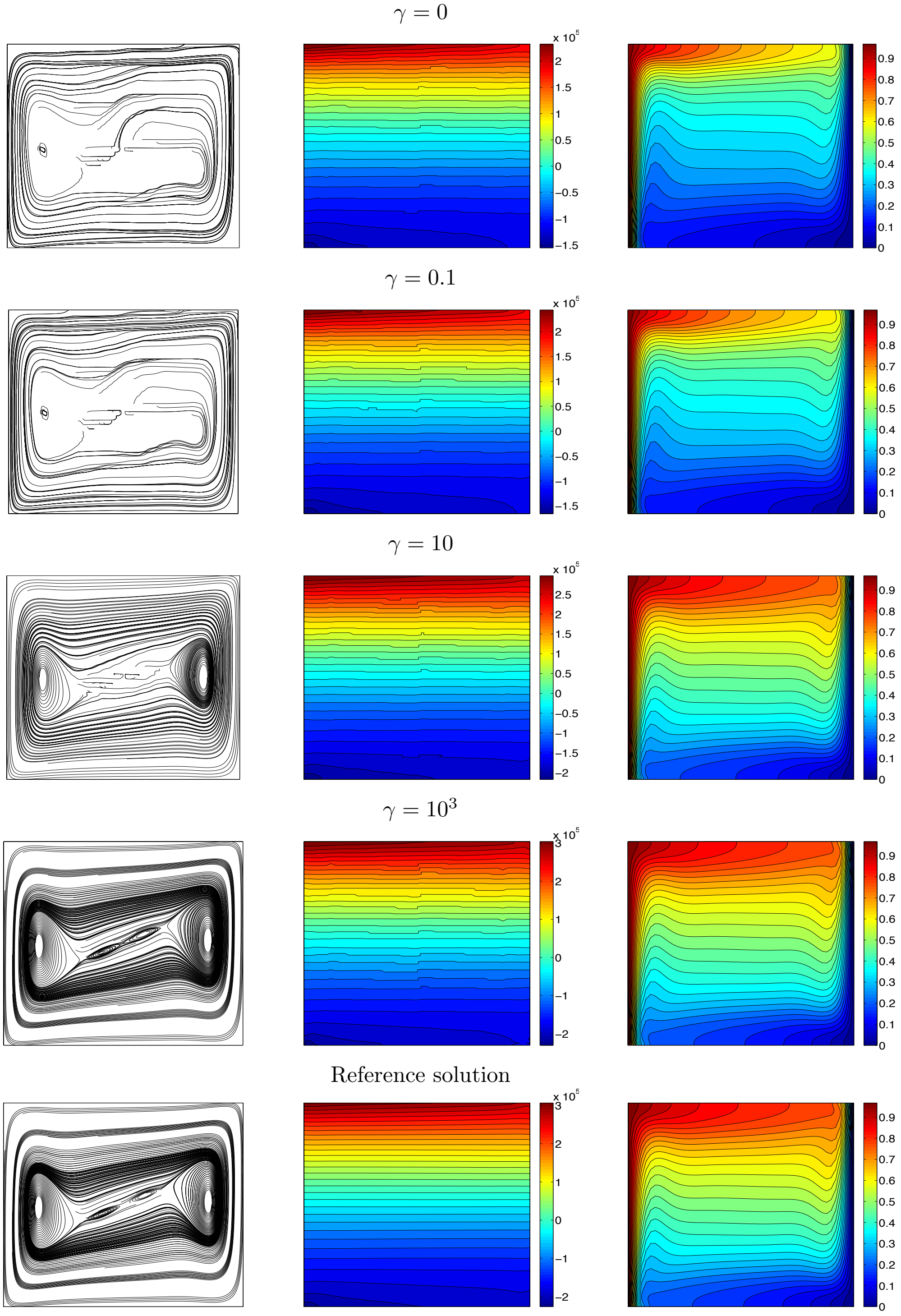

Figure 3: Shown above are plots of CR solutions' velocity streamlines (left), pressure contours (middle), and temperature contours(right), for varying $\gamma$, along with a reference solution at the bottom.

## Acknowledgments

The authors would especially like to thank Christoph Lehrenfeld for several related fruitful discussions on stabilization and hybridization and the invaluable help he provided in using the finite element library `NGSolve` in the context of this work.